\documentclass[twoside,11pt]{amsart}

\usepackage{amsmath,latexsym,amssymb, times, enumerate, hyperref}
\usepackage{graphicx,supertabular,paralist}
\usepackage[dvips,all]{xy}

\usepackage{tikz}
\usetikzlibrary{patterns}

\input amssym.def
\input amssym
\input xypic
\input xy
\xyoption{all}
\def\demo{\noindent{\bf Proof. }}
\def\sqr#1#2{{\vcenter{\hrule height.#2pt
        \hbox{\vrule width.#2pt height#1pt \kern#1pt
                \vrule width.#2pt}
        \hrule height.#2pt}}}
\def\square{\mathchoice\sqr64\sqr64\sqr{4}3\sqr{3}3}
\def\QED{\hfill$\square$}

\def\p{{\mathfrak p}}

\newtheorem{Theorem}{Theorem}[section]

\newtheorem{Lemma}[Theorem]{Lemma}
\newtheorem{Corollary}[Theorem]{Corollary}
\newtheorem{Proposition}[Theorem]{Proposition}
\newtheorem{Conjecture}[Theorem]{Conjecture}

\newtheorem{Discussion}[Theorem]{Discussion}
\newtheorem{Notation}[Theorem]{Notation}
\newtheorem{Assumptions and Discussion}[Theorem]{Assumptions and Discussion}

\newtheorem{Remark}[Theorem]{Remark}
\newtheorem{Example}[Theorem]{Example}
\newtheorem{Definition}[Theorem]{Definition}

\newtheorem{Algorithm}[Theorem]{Algorithm}

\newcommand{\pd}{\mathop{\mathrm{pd}}\nolimits}

\newcommand{\nr}{\mathop{\mathrm{nr}}\nolimits}
\renewcommand{\H}{\mathcal{H}}
\newcommand{\Q}{\mathcal{Q}}

\renewcommand{\SS}{\mathcal{S}}

\newcommand{\BB}{\mathcal{B}}

\newcommand{\G}{\mathcal{G}}
\newcommand{\HH}{\mathcal{H}}

\def\ovl#1{{\overline{#1}}}
\def\ww#1{\widetilde{#1}}

\begin{document}

\baselineskip=16pt

\baselineskip=16pt

\title[Hypergraphs with high projective dimension and 1-dimensional Hypergraphs]
{\Large\bf Hypergraphs with high projective dimension and 1-dimensional Hypergraphs}

\author[K.-N. Lin, and P. Mantero]
{K.-N. Lin and P. Mantero}

\thanks{AMS 2010 {\em Mathematics Subject Classification}.
Primary 13D02; Secondary 05E40.}


\address{
Penn State Greater Allegheny, Academic Affairs,  McKeesport, PA
}
\email{kul20@psu.edu}

\address{University of Arkansas, Department of Mathematical Sciences, Fayetteville, AR
} \email{pmantero@uark.edu}

\

\date{\today}

\vspace{-0.1in}

\begin{abstract}
We prove a sufficient and a necessary condition for a square-free monomial ideal $J$ associated to a (dual) hypergraph to have projective dimension equal to the minimal number of generators of $J$ minus 2. We also provide an effective explicit procedure to compute the projective dimension of 1-dimensional hypergraphs $\H$ when each connected component contains at most one cycle. An algorithm to compute the projective dimension is also included. Applications of these results are given; they include, for instance, computing the projective dimension of monomial ideals whose associated hypergraph has a spanning Ferrers graph.
\end{abstract}

\maketitle

\section{Introduction}
Let $k$ be a field, $R=k[x_{1},\cdots,x_{n}]$ a polynomial ring with indeterminates $x_{1},\ldots,x_{n}$ and let $I$ be a homogeneous ideal of $R$. Over the last decades there has been great interest in determining (or bounding) two fundamental invariants of $I$, the projective dimension $\pd(R/I)$ and  the Castelnuovo-Mumford regularity ${\rm reg}(J)$. These two invariants play an important role in algebraic geometry, commutative algebra and combinatorial algebra. 
To compute these two invariants, it is natural to determine the minimal graded free resolution of $I$ which, however, is often difficult and computationally expensive to find. A slightly different approach consists in finding upper bounds for these two invariants for $I$, by computing the projective dimension and the regularity of one of its initial ideals $J={\rm in}_{\tau}(I)$. The projective dimension and the regularity of a monomial ideal are preserved by polarization, thus it is sufficient to consider square-free monomial ideals.
In general, computing the regularity ${\rm reg}(J)$ can be hard and computationally very expensive; for square-free monomomial ideals, however, one can take advantage of the equality ${\rm reg}(J)=\pd(R/J^{\vee})$, 
where $J^{\vee}$ is the Alexander dual of $J$, and reduce the problem to computing the projective dimension of a square-free monomial ideal, which is then 
an active area of research.

In general, given a square-free monomial ideal $J$, several combinatorial structures can be associated to it (e.g. 
simplicial complexes, graphs, hypergraphs or dual hypergraphs). They have been consistently used to establish combinatorial characterizations for the projective dimension or regularity of $J$ under additional assumptions, see -- among the many papers on the subject -- \cite{C},\cite{DHS},\cite{DS},\cite{DS1},\cite{F},\cite{H},\cite{HW},\cite{Ku},\cite{LMc},\cite{MV}.
In the present paper, we employ the combinatorial structures of ``dual hypergraphs'' in the sense of \cite{Be} (which we call for simplicity hypergraphs) to determine $\pd(R/J)$ for classes of monomial ideals $J$ whose associated (dual) hypergraphs satisfy certain combinatorial assumptions. We recall that the association of a dual hypergraph to a monomial ideal $J$ was first introduced by Kimura, Terai and Yoshida, who employed it to compute the arithmetical rank of certain square-free monomial ideals \cite{KTY}. In the last few years, various work has been done to determine invariants or detect properties of $J$ using this combinatorial association, see for example \cite{HL},\cite{KM},\cite{KRT}, \cite{KTY1},\cite{LMa},\cite{LMc}. 

One of the main results of this paper is motivated and, in some sense, is the continuation of work of Kimura, Rinaldo and Terai, who found necessary and sufficient combinatorial conditions on the hypergraph of $J$ to have $\pd(R/J)=\mu(J)-1$, i.e.  projective dimension equal to its minimal number of generators of $J$ minus 1. In this paper we provide a sufficient condition for $\pd(R/J)=\mu(J)-2$; if, in addition, the hypergraph associated to $J$ is a bipartite graph, then a necessary condition is also given (Theorem \ref{pdH-2}).
As an application, we compute the projective dimension of any $J$ whose 1-dimensional subhypergraph has a spanning Ferrers graph (Corollary \ref{BabyF}).

The second main result of this paper is the continuation of authors' previous work \cite{LMa}, where a combinatorial formula for $\pd(R/J)$ was found when the hypergraph associated to $J$ is a string or a cycle. In the present paper, we find the projective dimension of $\pd(R/J)$ when its (dual) hypergraph is a disjoint
union of trees or graphs containing at most one cycle (Theorem \ref{uniqueCycle}). In the special case of a disjoint union of trees, also results of Morey and Villarreal, and Faridi  apply \cite{MV} \cite{F} (because our tree hypergraphs are simplicial trees and then $J$ is, in these cases, sequentially Cohen-Macaulay); their results state that $\pd(R/J)$ equals the {\it big height} of $J$, i.e. the largest height of an associated prime of $J$. In this scenario our combinatorial result provides an alternative way to the above-mentioned algebraic formula for $\pd(R/J)$. When an explicit irredundant primary decomposition of $J$ is given, the big height of $J$ is easily computed; on the other hand, when the combinatorial structure is given, our formula usually provides a faster way to compute $\pd(R/J)$, especially useful when $J$ involves a large number of variables (because the big height is computed as the maximum of all possible vertex covers of the corresponding simplicial structure).

The key idea for this result is to develop a process for breaking a ``large'' ideal into ``smaller'' ideals having disjoint combinatorial structures, thus reducing the computational cost of finding $\pd(R/J)$ (Propositions \ref{AttachString} and \ref{cutTree}).
Another consequence of Theorem \ref{uniqueCycle} is a combinatorial formula for $\pd(R/J)$ when the associated combinatorial structure can be described in terms of ``small'' stars; in these cases the formula has a flavor similar
to the main result of \cite{LMa} (Propositions \ref{EstarsAndStars} and \ref{StarString}).

The paper is organized as follows: in Section 2 we set the ground work for the paper, establish notations, review properties and prove a few additional tools employed in the later sections.
Section 3 is developed around the first main theorem, featuring the necessary and the sufficient condition
to have $\pd(R/J)=\mu(J)-2$; it also contains an application to hypergraphs with a spanning Ferrers graph. In Section 4 we introduce an argument which essentially allows us to replace a large 1-dimensional hypergraph with the disjoint union of smaller hypergraphs; we employ it to prove Theorem \ref{uniqueCycle} and provide a few applications. In Section 5 we have included an algorithm to compute the projective dimension of the connected hypergraphs to which one can apply Theorem \ref{uniqueCycle}.

\section{Background and a few lemmas}

We recall that the definition of (dual) hypergraph on the vertex set $V=[\mu]=\{1,2,\ldots,\mu\}$, see \cite{Be}.
\begin{Definition}\label{H}
A {\em (dual) hypergraph} on $V=[\mu]$  is a subset $\H$ of the power set $\mathcal P(V)$ such that $\bigcup\limits_{F\in \H}F=V$.
$\H$ is {\em separated} if, moreover, for every $1\leq j_1 <\, j_2\leq \mu$, there exist faces $F_1,F_2\in \H$ with $j_1\in F_1 \cap (V\setminus F_2)$ and $j_2\in F_2 \cap (V\setminus F_1)$.
\end{Definition}

Let $R=k[x_1,...,x_n]$ be a polynomial ring over a field $k$. If $I$ is a square-free monomial ideal in $R$, then one can associate a separated hypergraph $\H(I)$ to it: let $m_1,\ldots,m_{\mu}$ be the minimal monomial generating set for $I$, the hypergraph $\H(I)$ is defined as 
$$\H(I):=\{\{j \in V:x_i\,|\,m_j\}:i=1,2,\cdots ,n \}.$$

The hypergraph $\H(I)$ defined above is sometimes called the dual hypergraph of $I$ and should not be confused with the hypergraph constructed from $I$ by setting as vertices the variables of the polynomial ring, and having the faces correspond to the generators of the ideal. Also, following \cite{HL},\cite{KRT},\cite{KTY},\cite{LMa},\cite{LMc} we assume all the hypergraphs are separated, unless otherwise stated.

\begin{Definition}\label{not}
Let $I$ be a square-free monomial ideal with minimal monomial generating set $\{m_1,\cdots,m_{\mu}\}$. We set $\pd(\H)$ for $\pd(R/I)$, where $\H=\H(I)$ is the hypergraph associated to the square-free monomial ideal $I$, and call it the {\em projective dimension} of $\H$.
\end{Definition}

Conversely, given a separated hypergraph $\H$ with vertex set $V=[\mu]$, one can associate to it multiple monomial ideals, see for instance, \cite{KRT} or \cite{LMc}. In our proofs we will always associate to $\H$ a (standard) square-free monomial $I(\H)$ minimally generated by monomials $m_1 ,...m_{\mu}$ with the additional property that for every face $F$ in $\H$, there is a unique variable $x_F$ such that $x_F|m_l$ if and only if $l$ is in $F$.
This can be done without loss of generality, since in \cite[Proposition~2.2 and Corollary~2.4]{LMa}, the authors showed that any two square-free monomial ideals associated to the same separated hypergraph $\H$ have the same Betti numbers and projective dimension.

We now summarize a few combinatorial operations and their algebraic counterparts.

\begin{Definition}\label{deform}
Let $\H$ be a hypergraph, $I=I(\H)\subseteq R$ be the (standard) square-free monomial ideal associated to it, let $F$ be a face in $\H$ and $x_F\in R$ be the variable associated to $F$; also, let $v$ be a vertex in $\H$ and $m_v\in I$ be the monomial generator associated to it. We define the following operations on $\H$.
\begin{itemize}
\item[(i)] The hypergraph $\H_v$ obtained by {\rm remotion} of $v$ from $\H$ is defined as follows: let $A$ be the set obtained by removing $m_v$ from the set of minimal monomial generators of $I$, set $I_v=(m\,|\,m\in A)$ , then $\H_v=\H(I_v)$; iterating this operation, one writes $\H_{v_1,\ldots,v_r}$ for the hypergraph obtained by removing multiple vertices $v_1,\ldots,v_r$;

\item[(ii)] the hypergraph $\H_v:v=\Q_v$ is the hyergraph $\H(I_v:m_v)$ where $I_v$ and $\H_v$ are as in (i).

\item[(iii)] the hypergraph $\H:F$ obtained by {\rm cancellation} of $F$ in $\H$ is the hypergraph associated to $I:x_F$;

\item[(iv)] the hypergraph $\H^F$ obtained from $\H$ by {\rm cutting} $F=\{v_{i_1},\ldots,v_{i_r}\}$ is defined as follows: assume $m_{i_j}$ is the monomial in $I=I(\H)\subseteq R=k[y,x_F]$ corresponding to the vertex $v_{i_j}$; now set $R'=k[y,x_{F_1},\ldots,x_{F_r}]$ where $x_{F_1},\ldots,x_{F_r}$ are new variables, and consider the monomial ideal $I^F\subseteq R'$ obtained from $I=I(\H)$ by changing only the monomial generators $m_{i_j}$ as follows: replace $m_{i_j}$ by $m_{i_j}'=\frac{x_{F_j}}{x_F}m_{i_j}$. The hypergraph $\H^F$ is $\H(I^F)$.
\end{itemize} 
\end{Definition}

\begin{Example}
In Figure \ref{deformH} we fix the hypergraph $\mathcal{H}$, a vertex $v$, and faces $F$ and $E$. The hypergraphs $\mathcal{H}_{v}$, $\mathcal{H}_{v}:v=\mathcal{Q}_v$,\; $\mathcal{H}:F$, and $\mathcal{H}^{E}$ are represented in Figure \ref{DeformHv}.\end{Example}

\begin{figure}[h] 
\caption{}\label{deformH}
\begin{center}
\begin{tikzpicture}[thick, scale=1.3]

\draw  [shape=circle] (0,0) circle (.1);
\draw  [shape=circle] (1,0) circle (.1);
\draw  [shape=circle] (2,0) circle (.1);
\draw  [shape=circle] (3,0) circle (.1);
\shade [shading=ball, ball color=black]  (4,-0.5) circle (.1);
\draw  [shape=circle] (0,-1) circle (.1);
\draw  [shape=circle] (1,-1) circle (.1);
\draw  [shape=circle] (2,-1) circle (.1);
\shade [shading=ball, ball color=black] (0.5,0.5) circle (.1) node [above] {$v$};

\draw [line width=1.2pt ] (0.5,0.5)--(1,0)  ;

\draw [line width=1.2pt ] (0,0)--(0,-1)  ;
\draw [line width=1.2pt ] (0,0)--(1,-1)  ;
\draw [line width=1.2pt ] (0,0)--(2,-1)  ;

\draw [line width=1.2pt ] (1,0)--(0,-1)  ;
\draw [line width=1.2pt ] (1,0)--(1,-1)  ;
\draw [line width=1.2pt ] (2,0)--(0,-1)  ;
\draw [line width=1.2pt ] (2,0)--(1,-1)  ;
\draw [line width=1.2pt ] (3,0)--(0,-1)  ;
\draw [line width=1.2pt ] (3,0)--(4,-0.5)  ;
\draw [line width=1.2pt ] (4,-0.5)--(2,-1)  ;

\path [pattern=north east lines, pattern color=red]     (2,0)--(3,0)--(2,-1)--(1,-1)--cycle;
\path [pattern=north east lines, pattern color=blue]     (0.5,0.5)--(0,0)--(1,0)--cycle;

\path (3.5,0.5)--(3.5,-0.5) node [pos=.5, right ] {$F$};
\path (2,0)--(2,-1) node [pos=.5, right ] {$E$};

\end{tikzpicture}
\end{center}
\end{figure}

\begin{figure}[h] 
\caption{}\label{DeformHv}
\begin{center}
\begin{tikzpicture}[thick, scale=1.3]

\draw  [shape=circle] (0,0) circle (.1);
\shade [shading=ball, ball color=black] (1,0) circle (.1);
\draw  [shape=circle] (2,0) circle (.1);
\draw  [shape=circle] (3,0) circle (.1);
\shade [shading=ball, ball color=black]  (4,-0.5) circle (.1);
\draw  [shape=circle] (0,-1) circle (.1);
\draw  [shape=circle] (1,-1) circle (.1);
\draw  [shape=circle] (2,-1) circle (.1);

\draw [line width=1.2pt ] (0,0)--(0,-1)  ;
\draw [line width=1.2pt ] (0,0)--(1,-1)  ;
\draw [line width=1.2pt ] (0,0)--(2,-1)  ;
\draw [line width=1.2pt ] (0,0)--(1,0)  ;

\draw [line width=1.2pt ] (1,0)--(0,-1)  ;
\draw [line width=1.2pt ] (1,0)--(1,-1)  ;
\draw [line width=1.2pt ] (2,0)--(0,-1)  ;
\draw [line width=1.2pt ] (2,0)--(1,-1)  ;
\draw [line width=1.2pt ] (3,0)--(0,-1)  ;
\draw [line width=1.2pt ] (3,0)--(4,-0.5)  ;
\draw [line width=1.2pt ] (4,-0.5)--(2,-1)  ;

\path [pattern=north east lines, pattern color=red]     (2,0)--(3,0)--(2,-1)--(1,-1)--cycle;

\path (3.5,0.5)--(3.5,-0.5) node [pos=.5, right ] {$F$};
\path (2,0)--(2,-1) node [pos=.5, right ] {$E$};

\path (1.5,-1.5)--(1.5,-1.5) node [pos=.5, right ] {$\mathcal{H}_v$};

\draw  [shape=circle] (5,0) circle (.1);
\draw  [shape=circle] (6,0) circle (.1);
\draw  [shape=circle] (7,0) circle (.1);
\draw  [shape=circle] (8,0) circle (.1);
\shade [shading=ball, ball color=black]  (9,-0.5) circle (.1);
\draw  [shape=circle] (5,-1) circle (.1);
\draw  [shape=circle] (6,-1) circle (.1);
\draw  [shape=circle] (7,-1) circle (.1);

\draw [line width=1.2pt ] (5,0)--(5,-1)  ;
\draw [line width=1.2pt ] (5,0)--(6,-1)  ;
\draw [line width=1.2pt ] (5,0)--(7,-1)  ;

\draw [line width=1.2pt ] (6,0)--(5,-1)  ;
\draw [line width=1.2pt ] (6,0)--(6,-1)  ;
\draw [line width=1.2pt ] (7,0)--(5,-1)  ;
\draw [line width=1.2pt ] (7,0)--(6,-1)  ;
\draw [line width=1.2pt ] (8,0)--(5,-1)  ;
\draw [line width=1.2pt ] (8,0)--(9,-0.5)  ;
\draw [line width=1.2pt ] (9,-0.5)--(7,-1)  ;

\path [pattern=north east lines, pattern color=red]     (7,0)--(8,0)--(7,-1)--(6,-1)--cycle;

\path (8.5,0.5)--(8.5,-0.5) node [pos=.5, right ] {$F$};
\path (7,0)--(7,-1) node [pos=.5, right ] {$E$};

\path (6.5,-1.5)--(6.5,-1.5) node [pos=.5, right ] {$\mathcal{Q}_v$};

\draw  [shape=circle] (0,-2.5) circle (.1);
\draw  [shape=circle] (1,-2.5) circle (.1);
\draw  [shape=circle] (2,-2.5) circle (.1);
\draw  [shape=circle] (3,-2.5) circle (.1);
\shade [shading=ball, ball color=black]  (4,-3) circle (.1);
\draw  [shape=circle] (0,-3.5) circle (.1);
\draw  [shape=circle] (1,-3.5) circle (.1);
\draw  [shape=circle] (2,-3.5) circle (.1);
\shade [shading=ball, ball color=black] (0.5,-2) circle (.1) node [above] {$v$};

\draw [line width=1.2pt ] (0.5,-2)--(1,-2.5)  ;

\draw [line width=1.2pt ] (0,-2.5)--(0,-3.5)  ;
\draw [line width=1.2pt ] (0,-2.5)--(1,-3.5)  ;
\draw [line width=1.2pt ] (0,-2.5)--(2,-3.5)  ;

\draw [line width=1.2pt ] (1,-2.5)--(0,-3.5)  ;
\draw [line width=1.2pt ] (1,-2.5)--(1,-3.5)  ;
\draw [line width=1.2pt ] (2,-2.5)--(0,-3.5)  ;
\draw [line width=1.2pt ] (2,-2.5)--(1,-3.5)  ;
\draw [line width=1.2pt ] (3,-2.5)--(0,-3.5)  ;

\draw [line width=1.2pt ] (4,-3)--(2,-3.5)  ;

\path [pattern=north east lines, pattern color=red]     (2,-2.5)--(3,-2.5)--(2,-3.5)--(1,-3.5)--cycle;
\path [pattern=north east lines, pattern color=blue]     (0.5,-2)--(0,-2.5)--(1,-2.5)--cycle;

\path (2,-2.5)--(2,-3.5) node [pos=.5, right ] {$E$};

\path (1.5,-4)--(1.5,-4) node [pos=.5, right ] {$\mathcal{H}:F$};

\draw  [shape=circle] (5,-2.5) circle (.1);
\draw  [shape=circle] (6,-2.5) circle (.1);
\shade [shading=ball, ball color=black] (7,-2.5) circle (.1);
\shade [shading=ball, ball color=black] (8,-2.5) circle (.1);
\shade [shading=ball, ball color=black]  (9,-3) circle (.1);
\draw  [shape=circle] (5,-3.5) circle (.1);
\shade [shading=ball, ball color=black](6,-3.5) circle (.1);
\shade [shading=ball, ball color=black] (7,-3.5) circle (.1);
\shade [shading=ball, ball color=black] (5.5,-2) circle (.1) node [above] {$v$};

\draw [line width=1.2pt ] (5.5,-2)--(6,-2.5)  ;
\draw [line width=1.2pt ] (5,-2.5)--(5,-3.5)  ;
\draw [line width=1.2pt ] (5,-2.5)--(6,-3.5)  ;
\draw [line width=1.2pt ] (5,-2.5)--(7,-3.5)  ;

\draw [line width=1.2pt ] (5,-3.5)--(6,-2.5)  ;
\draw [line width=1.2pt ] (5,-3.5)--(7,-2.5)  ;
\draw [line width=1.2pt ] (5,-3.5)--(8,-2.5)  ;
\draw [line width=1.2pt ] (6,-2.5)--(6,-3.5)  ;
\draw [line width=1.2pt ] (8,-2.5)--(9,-3)  ;
\draw [line width=1.2pt ] (6,-3.5)--(7,-2.5)  ;

\path [pattern=north east lines, pattern color=blue]     (5.5,-2)--(5,-2.5)--(6,-2.5)--cycle;

\path (8,-2.5)--(9,-2.5) node [pos=.5, right ] {$F$};

\path (6,-4)--(7,-4) node [pos=.5, right ] {$\mathcal{H}^E$};

\end{tikzpicture}
\end{center}
\end{figure}

\begin{Discussion}
We now discuss and explain briefly the operations defined in Definition \ref{deform}.
\begin{itemize}
\item[(i)] the ``remotion'' of $v$ corresponds, in the realm of simplicial complexes, to taking the sub-simplicial complex obtained by removing the face associated to $v$. Here we call it remotion because from the point of view of (dual) hypergraphs it corresponds to removing $v$ from $\H$ and contracting the faces containing $v$.
\item[(ii)] the operation of ``cancelling'' $F$ corresponds, in the realm of simplicial complexes, to the operation of contraction of the vertex associated to $x_F$. Here we call it cancelling because from the point of view of (dual) hypergraphs it corresponds to cancelling $F$ from $\H$.
\item[(iii)] The operation of ``cutting'' derives its name from its combinatorial meaning, because $\H^F$ can be interpreted as the hypergraph obtained by cutting the face $F$ into $r$ parts (one for each vertex of $F$) and retracting each of them back to the corresponding vertex $v_{i_j}$. This is different from the cancellation of $F$ as long as one the vertices $v_{i_j}$ is open in $\H$, because after this operation $v_{i_j}$ will become closed.
\end{itemize}
\end{Discussion}

We now relate $\pd(\H)$ with the projective dimension of the hypergraphs defined in Definition \ref{deform}. 
\begin{Lemma}\label{ref}\cite[Lemma~2.6]{LMa}
Let $\H$ be a hypergraph. If $\{v\}\in \H$, then $$\pd(\H)={\rm max}\{\pd(\H_v),\pd(\Q_v)+1\}.$$
\end{Lemma}

\begin{Lemma}\label{red}\cite[Lemma~2.11]{LMa}
Let $\H$ be a hypergraph. If $\{v\}\in \H$ and all its neighbors are closed vertices, then $\pd(\H)=\pd(\H_v)+1$.
\end{Lemma}

\begin{Lemma}\label{loc}\cite[Lemma~2.8]{LMa}
If $\H'\subseteq\H$ are hypergraphs with $\mu(\H')=\mu(\H)$, then $\pd(\H')\leq\pd(\H)$.
\end{Lemma}

\begin{Proposition}\label{eq}\cite[Proposition~2.10]{LMa}
Let $\H',\H$ be hypergraphs with $\H=\H'\cup F$ where $F=\{i_1,\ldots,i_r\}$. If $\{i_j\}\in \H'$ for all $j$, then $\pd(\H')=\pd(\H:F)=\pd(\H)$.
\end{Proposition}

We recall the following folklore fact that can be proved, for instance, by means of Taylor's resolution \cite{T}.
\begin{Remark}\label{J+1}
Let $\H$ be a hypergraph then $\pd(\H_v)\leq\pd(\H)\leq\pd(\H_v)+1$.
\end{Remark}

In the following, we will need to know how can the projective dimension of a hypergraph vary if we make an open vertex become closed. This is studied in the following results.

\begin{Lemma}\label{H^0}
Let $\H$ be a hypergraph and let $\H^0$ be the hypergraph obtained by making one closed vertex $v$ in $\H$ become open.
Then $\pd(\H^0)\leq \pd(\H)\leq \pd(\H^0)+1$.
\end{Lemma}

\demo
The inequality on the left follows by Lemma \ref{loc}. Notice that $\H=\H^0\cup \{v\}$ and $\H_v=\H^0_v$. Then, by Remark \ref{J+1}, we have $\pd(\H)\leq\pd(\H_v)+1=\pd(\H^0_v)+1$. The desired inequality now follows because $\pd(\H^0_v)\leq \pd(\H^0)$, by Remark \ref{J+1}.
\QED 
\bigskip

We note that the hypergraph $\H:F$ obtained by cancelling $F$ can also be obtained by localization.
\begin{Lemma}\label{local}
Let $F$ be a face of a hypergraph $\H$, let $I=I(\H)\subseteq R$, let $\p$ be the ideal of all variables in $R$ except $x_F$, let $\widetilde{R}=R_{x_F}$, $\widetilde{I}=I\widetilde{R}$, $S={\rm gr}_{\p\widetilde{R}}(\widetilde{R})$, and let  $I_1$ be the ideal of initial forms of $\widetilde{I}$ in $S$. 

If $\H_1=\H(I_1)$, then $\H:F=\H_1$. In particular $\pd(\H:F)\leq \pd(\H)$.
\end{Lemma}

\demo
Let $R=k[y_1,\ldots,y_n,x_F]$. Since $I$ is square-free, we can write $I=x_FJ+K$ for some square-free monomial ideals $J,K$ of the form $J=J'T$ and $K=K'T$ with $T=k[y_1,\ldots,y_n]$. In particular, $x_F$ is regular on $R/K$. Then, it is easily checked that the ideal associated to $\H:F$ is
$$I:x_F=(x_FJ+K):x_F=J+K.$$
With the same notation as above, we have $\widetilde{I}=(J+K)\widetilde{R}$, therefore $I_1=(J+K)S$, hence the hypergraph associated to $I_1$ is the same as the hypergraph associated to $J+K$.\\
In particular, we obtain $\pd(R/I)\geq \pd(R_{\p}/I_{\p})=\pd(S/I_1)=\pd(R/I:x_F)$.
\QED 
\bigskip

\begin{Lemma}\label{H^F}
Let $\H$ be a hypergraph and $F$ be a face of $\H$ then $\pd(\H)\leq\pd(\H^F)$.
\end{Lemma}

\demo
Let $\H'=\H^F\cup F$ then $\pd(\H')=\pd(\H^F)$ by the definition of $\H^F$ and Proposition \ref{eq}. The conclusion now follows by Lemma \ref{loc} since  $\mu(\H)=\mu(\H')$, and $\H\subseteq\H'$.
\QED 
\bigskip

Also, we obtain a possibly useful criterion to compute $\pd(\H)$ when one more piece of information is known. 
\begin{Corollary}\label{different}
Let $\H$ be a hypergraph and let $\H^0$ be the hypergraph obtained by making one closed vertex $v$ become open. If $\pd(\H^0)\neq \pd(\H_v)$, then $$\pd(\H)={\rm max}\{\pd(\H^0),\pd(\H_v)\}={\rm min}\{\pd(\H^0),\pd(\H_v)\}+1.$$
\end{Corollary}
\demo
If $\pd(\H^0)<\,\pd(\H_v)$, then by Remark \ref{J+1} and Lemma \ref{H^0} we have
$$\pd(\H^0)+1\leq \pd(\H_v)\leq \pd(\H)\leq \pd(\H^0)+1$$
which yields $\pd(\H)=\pd(\H_v)=\pd(\H^0)+1$. The case where $\pd(\H_0)>\pd(\H_v)$ is proved symmetrically.
\QED
 \bigskip

\section{Large projective dimension and Generalized Ferrer Graphs}

Let us recall that two vertices $v\neq w$ in $\H$ are {\it neighbors} if there is a face of $\H$ containing both of them. Let nb$(v)$ denote the set of all neighbors of the vertex $v$, its cardinality $\deg(v)$ is called the {\em degree} of $v$. If $\deg(v)=0$, i.e. $v$ has no neighbors, then $v$ is called {\it isolated}; in this case one has $\pd(\H)=\pd(\H')+1$, where $\H'=\H\setminus \{v\}$. Thus, each isolated vertex contributes to the projective dimension with one unit. Since our focus is on the projective dimension of the hypergraphs, and the projective dimension of a hypergraph with two disconnected subhypergraphs is the sum of the projective dimensions of the subhypergraphs, we may assume all hypergraphs have no isolated vertices.

\begin{Notation} 
Let $\H$ be a hypergraph with vertex set $V(\H)$, following \cite{KRT} we write
\begin{itemize}
\item $W(\H)=\{i\in V|\{i\}\notin \H\}$ for the open vertex set of $\H$,
\item $\H_{U}=\{F\in \H\,:\, F\subseteq U\}$ for the restriction to a subset $U\subseteq V(\H)$ of the vertex set;
\item $\H^i=\{F\in \H: \dim F\leq i\}$ for the $i$-th dimensional subhypergraph of $\H$.
\end{itemize}
Moreover, for a subset $U\subseteq V(\H)$ we define $\H_{\overline{U}}=\H\cup \{ \{i\}\,:\,i \notin U\}$ for the hypergraph obtained by making all vertices of $\H$ not in $U$ become closed.
\end{Notation}

It is well-known that if $W(\H)=\emptyset$, then $\pd(\H)=|V(\H)|$. The following theorem by Kimura, Rinaldo and Terai characterizes when $\pd(\H)=|V(\H)|-1$.

\begin{Theorem}\cite[Theorem~4.3]{KRT}\label{KRT}
Let $\H$ be a hypergraph, then $\pd(\H)=|V(\H)|-1$ if and only if $\H$ satisfies the following condition
$$\begin{array}{ll}
(\star) & W(\H)\neq\emptyset \mbox{ and either the 1-dimensional part } \H_{W(\H)}^1 \mbox{ of }\H_{W(\H)} \mbox{ contains a spanning complete}\\
& \mbox{ bipartite graph, or there is a vertex }v \mbox{ such that } \{\{v,w\}\in \H\,|\,\mbox{ for every } w\in W(\H)\}.
\end{array}$$ 
\end{Theorem}

One then has the following corollary.

\begin{Corollary}\label{CKRT}
Let $\H$ be a hypergraph. Then $\pd(R/\H)\leq |V(\H)|-2$ if and only if $\H$ satisfies the following condition
$$\begin{array}{ll}
(\star\star) & W(\H)\neq\emptyset,\; \H_{W(\H)}^1 \mbox{ does not contain a spanning complete bipartite graph, and}\\
& \mbox{ there is no vertex }v \mbox{ such that } \{\{v,w\}\in \H\,|\,\mbox{ for every } w\in W(\H)\}.
\end{array}$$ \end{Corollary}

Then, the next step is trying to determine the hypergraphs with $\pd(\H)=|V(\H)|-2$. We define the following assumption:\\
\\
\noindent\begin{tabular}{ll}
$(\sharp)$ & $\H$ satisfies $(\star\star)$ and there is a partition $\{V_1,V_2\}$ of the vertex set $V(\H)$ such that both\\ & $\H_{\overline{V_1}}$ and $\H_{\overline{V_2}}$ satisfy $(\star)$. 
\end{tabular}
\\

We now prove that $(\sharp)$ gives, in general, a sufficient condition for $\pd(\H)=|V(\H)|-2$. For 1-dimensional bipartite hypergraphs $\H$ we prove a necessary condition for $\pd(\H)=|V(\H)|-2$, which is very similar to $(\sharp)$.
\begin{Theorem}\label{pdH-2}
Let $\H$ be a hypergraph.\\
$($i$)$ If $\H$ satisfies $(\sharp)$, then $\pd(\H)=|V(\H)|-2$.\\
$($ii$)$ If, furthermore, $\mathcal{H}$ is a 1-dimensional bipartite graph.
Then $\pd(\H)=V(\mathcal{H})-2$ implies there is a partition $\{V_1,V_2\}$ of the vertex set $V(\H)$ such that if $\mathcal{G}_{1}$ and $\mathcal{G}_{2}$ are obtained by cutting all edges of $\mathcal{H}$ between $V_{1}$
and $V_{2}$, then both $\mathcal{G}_{1}$ and $\mathcal{G}_{2}$ satisfy $(\star)$.
\end{Theorem}

\demo
(i) The inequality $\pd(\H)\leq|V(\H)|-2$ follows by assumption $(\sharp)$. To prove the other inequality, we observe that, after cancelling all the faces containing vertices both from $V_1$ and $V_2$, we are left with two disconnected subgraphs $\G_1=\H_{V_1}$ and $\G_2=\H_{V_2}$. Then, by Lemma \ref{local} we have 
$$\pd(\H)\geq \pd(\G_1)+\pd(\G_2)$$
We now show that, regardless of whether $\G_1$ and $\G_2$ are separated, one has $\pd(\G_i)=|V_i|-1$ for $i=1,2$. By symmetry, we only prove that $\pd(\G_1)=|V_1|-1$. Since $\H_{\overline{V_1}}$ satisfies property $(\star)$ and all open vertices of $\H_{\ovl{V_1}}$ are in $V_1$, if $\G_1\subseteq \H_{\overline{V_1}}$ is separated, then also $\G_1$ satisfies $(\star)$ and by Theorem \ref{KRT} one has $\pd(\G_1)=|V_1|-1$. 

We may then assume that $\G_1=\H_{V_1}$ is not separated. If $\H_{W(\H_{V_1})}^1$ contains a spanning complete bipartite graph where each bipartite set has more than one vertex, then $\G_1=\H_{V_1}$ is separated, which contradicts our assumption. We may then assume $\H_{\ovl{V_1}}$ has one vertex $w$ connected to all of its open vertices. Also, since $\G_1=\H_{V_1}$ is not separated, there are open vertices $v_1,...,v_t$ in $V_1$ which are vertices of faces containing vertices of both $\G_1$ and $\G_2$, and have the property that in $\H$ each $v_i$ has only one neighbor in the vertex set $V_1$. By the above, this neighbor is $w$ for every $i=1,\ldots,t$, i.e. they all have $w$ as a common neighbor.

Then, when we cancel the faces connecting vertices of $V_1$ with vertices of $V_2$, the hypergraph $\G_1$ just consists of $|V_1|-1$ closed vertices, because $w$ degenerates after the cancellation and all its neighbors, which include all open vertices of $\H_{V_1}$, become closed. Then $\G_1$ is saturated and thus $\pd(\G_1)=|V_1|-1$, whence the conclusion follows.

(ii) Let $W_{1}$ and $W_{2}$ be the two vertex sets of the bipartite
graph $\mathcal{H}$. Since $\mathcal{H}$ satisfies $(\star\star)$, then by Corollary \ref{CKRT}, there
are two open vertices, $v\in W_{1}$ and $w\in W_{2}$, which are not neighbors. Let $V_{1}=\{v,W_{2}\backslash w\}$
and $V_{2}=\{w,W_{1}\backslash v\}$, then we have $V_{1}\cup V_{2}=V(\mathcal{H})$,
and $V_{1}\cap V_{2}=\emptyset$. Since $\mbox{nb}(w)\cap V_{1}=\emptyset$
and $\mbox{nb}(v)\cap V_{2}=\emptyset$, then $v$ and $w$ are open after cutting
all edges of $\mathcal{H}$ between $V_{1}$ and $V_{2}$. Therefore
$\mbox{pd}(\mathcal{G}_{1})\leq|V_{1}|-1$ and $\mbox{pd}(\mathcal{G}_{2})\leq|V_{2}|-1$.
By Lemma \ref{H^F}, $\mbox{pd}(\mathcal{G}_{1})+\mbox{pd}(\mathcal{G}_{2})\geq\mbox{pd}(\mathcal{H})=|V(\mathcal{H})|-2=|V_{1}|+|V_{2}|-2$.
We conclude that $\mbox{pd}(\mathcal{G}_{1})=|V_{1}|-1$ and $\mbox{pd}(\mathcal{G}_{2})=|V_{2}|-1$,
and both $\mathcal{G}_{1}$ and $\mathcal{G}_{2}$ satisfies $(\star)$ 
by Theorem \ref{KRT}.  
\QED 
\bigskip

\begin{Remark}\label{v1v2}
It is easily seen that Theorem \ref{pdH-2}.(i) is also true (and has a much shorter proof) if one replaces $(\sharp)$ by the assumption that the restrictions $\H|_{V_1}$ and $\H|_{V_2}$ both satisfy $(\star)$ {\em and} are both separated. However, the requirement that $\H|_{V_1}$ and $\H|_{V_2}$ are separated is somewhat restrictive and, as we have proved above, unnecessary.

Also, the assumption that $\H_{\overline{V_1}}$ and $\H_{\overline{V_2}}$ have property $(\star)$ is much weaker than requiring that  $\H|_{V_1}$ and $\H|_{V_2}$ are separated sub-hypergraphs of $\H$ satisfying $(\star)$ as can be seen in a number of (even simple) examples. Consider, for instance, the 6-cycle graph $$\H=\{\{1\},\{1,2\},\{2,3\},\{3,4\},\{4\},\{4,5\},\{5,6\},\{6,1\}\}.$$
Then the vertex subsets $V_1=\{1,2,3\}$ and $V_2=\{4,5,6\}$  satisfy the assumptions of $(\sharp)$, thus $\pd(\H)=|V(\H)|-2=4$ by Theorem \ref{pdH-2}.(i); however, for every partition $\{U_1, U_2\}$ of $V(\H)$ neither $\H|_{U_1}$ nor $\H|_{U_2}$ is separated.
\end{Remark}

We suspect that the converse of Theorem \ref{pdH-2} (i) holds true provided $\H$ is 1-dimensional bipartite, although it does not follow by part (ii), because it is relatively easy to construct 1-dimensional hypergraphs $\H$ where a partition of $V(\H)$ constructed as in (ii) does not satisfy the assumption $(\sharp)$. However, in all the examples considered by the authors, we could always find another partition of $V(\H)$ satisfying $(\sharp)$. We then ask whether the following potential combinatorial characterization of 1-dimensional bipartite hypergraphs of projective dimension $|V(\H)|-2$ actually holds true:\\
\begin{Conjecture}\label{question}
Let $\H$ be a 1-dimensional bipartite hypergraph with $\pd(\H)\leq |V(\H)|-2$; then $\pd(\H)=|V(\H)|=2$ if and only if  $\H$ satisfies $(\sharp)$.
\end{Conjecture}

We have posed Conjecture \ref{question} under the additional assumption that $\H$ is a bipartite graph for two reasons: first, for the converse of Theorem \ref{pdH-2}.(i) one {\it needs} additional assumptions, as we show in Example \ref{ex}; and, secondly, because Theorem \ref{pdH-2}.(ii) shows that for bipartite graphs a condition very similar to $(\sharp)$ is indeed necessary.

\begin{Example}\label{ex}
The converse of Theorem \ref{pdH-2}.(i) does not hold in general, not even for graphs (i.e. 1-dimensional hypergraphs). For instance, let $\H$ be a 7-cycle graph whose vertices are all open. Then by the main result of \cite{LMa}, the projective dimension of $\H$ is $5=|V(\H)|-2$. However, one cannot find a partition $\{V_1,V_2\}$ of $V(\H)$ such that both $\H_{\overline{V_1}}$ and $\H_{\overline{V_1}}$ satisfy $(\star)$.
\end{Example}

As an application of Theorem \ref{pdH-2}, we show that  $\pd(\H)=|V(\H)|-2$ for a hypergraph $\H$ (not necessarily 1-dimensional) provided $\H_{W(\H)}^1$ has a spanning generalized Ferrers graph.
First, however, we recall the definition of generalized Ferrers graph. 
\begin{Definition}
A 1-dimensional bipartite graph $\{v_1,\ldots,v_s,w_1,\ldots,w_{\lambda_1}\}$, is a {\em  generalized Ferrers graph} if, after a permutation of vertices, there are two sequences of integers $\lambda=(\lambda_1,\ldots,\lambda_s)$ and $\tau=(\tau_{1},...,\tau_{t})$ such that
\begin{itemize}
\item $\lambda_{1}\geq\lambda_{2}\geq\cdots\geq\lambda_{s}>0$, 
\item $0=\tau_{1}\leq\tau_{2}\leq\cdots\leq\tau_{s}<\lambda_{s}$, 
\item $\lambda_i\leq \tau_i$ for every $i$,
\item and for every $i$, the vertex $v_i$ is connected to $w_{\tau_i+1},w_{\tau_i+2},\ldots,w_{\lambda_i}$ (in particular, $\lambda_i-\tau_i$ is the degree of $v_i$).
\end{itemize}
\end{Definition}

We give an example illustrating this definition.

\begin{Example}
Let $\H$ be a hypergraph with all open vertices whose 1-skeleton is described in Figure \ref{F}. Then, $\lambda=(7,7,6,5,4)$, $\tau=(0,0,1,1,2)$.

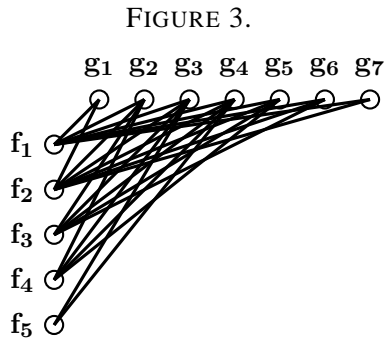
\begin{figure}[h] 
\caption{}\label{F}
\begin{center} 
\begin{tikzpicture} [thick, scale=0.6]

\draw  [shape=circle] (1,0) circle (.2) node [left] {$\bf{f_1}$};
\draw  [shape=circle] (1,-1) circle (.2) node [left] {$\bf{f_2}$};
\draw  [shape=circle] (1,-2) circle (.2) node [left] {$\bf{f_3}$};
\draw  [shape=circle] (1,-3) circle (.2) node [left] {$\bf{f_4}$};
\draw  [shape=circle] (1,-4) circle (.2) node [left] {$\bf{f_5}$};
\draw  [shape=circle] (2,1) circle (.2) node [above] {$\bf{g_1}$};
\draw  [shape=circle] (3,1) circle (.2) node [above] {$\bf{g_2}$};
\draw  [shape=circle] (4,1) circle (.2) node [above] {$\bf{g_3}$};
\draw  [shape=circle] (5,1) circle (.2) node [above] {$\bf{g_4}$};
\draw  [shape=circle] (6,1) circle (.2) node [above] {$\bf{g_5}$};
\draw  [shape=circle] (7,1) circle (.2) node [above] {$\bf{g_6}$};
\draw  [shape=circle] (8,1) circle (.2) node [above] {$\bf{g_7}$};

\draw [line width=1.2pt  ] (1,0)--(2,1)  ;
\draw [line width=1.2pt  ] (1,0)--(3,1)  ;
\draw [line width=1.2pt  ] (1,0)--(4,1)  ;
\draw [line width=1.2pt  ] (1,0)--(5,1)  ;
\draw [line width=1.2pt  ] (1,0)--(6,1)  ;
\draw [line width=1.2pt  ] (1,0)--(7,1)  ;
\draw [line width=1.2pt  ] (1,0)--(8,1)  ;
\draw [line width=1.2pt  ] (1,-1)--(2,1)  ;
\draw [line width=1.2pt  ] (1,-1)--(3,1)  ;
\draw [line width=1.2pt  ] (1,-1)--(4,1)  ;
\draw [line width=1.2pt  ] (1,-1)--(5,1)  ;
\draw [line width=1.2pt  ] (1,-1)--(6,1)  ;
\draw [line width=1.2pt  ] (1,-1)--(7,1)  ;
\draw [line width=1.2pt  ] (1,-1)--(8,1)  ;
\draw [line width=1.2pt  ] (1,-2)--(3,1)  ;
\draw [line width=1.2pt  ] (1,-2)--(4,1)  ;
\draw [line width=1.2pt  ] (1,-2)--(5,1)  ;
\draw [line width=1.2pt  ] (1,-2)--(6,1)  ;
\draw [line width=1.2pt  ] (1,-2)--(7,1)  ;
\draw [line width=1.2pt  ] (1,-3)--(3,1)  ;
\draw [line width=1.2pt  ] (1,-3)--(4,1)  ;
\draw [line width=1.2pt  ] (1,-3)--(5,1)  ;
\draw [line width=1.2pt  ] (1,-3)--(6,1)  ;
\draw [line width=1.2pt  ] (1,-4)--(4,1)  ;
\draw [line width=1.2pt  ] (1,-4)--(5,1)  ;

\end{tikzpicture}
\end{center}
\end{figure}
\end{Example}

\begin{Corollary}\label{BabyF}
Let $\H$ be a hypergraph satisfying $(\star\star)$. If $\H_{W(\H)}^1$ has a spanning generalized
Ferrers subgraph then $\pd(\H)=|V(\H)|-2$.
\end{Corollary}

\demo
After possibly a vertex permutation we may assume $\tau_1=\tau_2=\cdots=\tau_s=0$. Write $V(\H)=\{f_1,\ldots,f_s,g_1,...,g_{\lambda_1},z_1,\ldots,z_t\}$ where $\{f_1,..,f_s\}$ and $\{g_1,...,g_{\lambda_1}\}$ are the vertices corresponding to the two sets of generators of the bipartite graph spanning $\H_{W(\H)}^1$, and let $\{z_1,\ldots,z_t\}$ are all the closed vertices of $\H$. Take $V_1=\{\{g_1\}\cup \{f_i,|\, 2\leq i \leq s\}\}\cup \{z_k\,|\,1\leq k\leq t\}$, and $V_2=\{\{f_1\}\cup \{g_j\,|\, 2\leq j\leq \lambda_1\}\}$. Then $\H_{\overline{V_1}}$ and $\H_{\overline{V_2}}$ satisfy $(\star)$, thus the conclusion follows by Theorem \ref{pdH-2}.(i). 
\QED 
\bigskip

\begin{Example}
$\G_{1}$ is the green subgraph and $\G_{2}$ is the red subgraph in Figures \ref{FL}. Notice that $\H_{W(\H)}^1$ only needs to have a {\it  spanning} generalized Ferrers subgraph, so $\H$ could also contain higher dimensional faces and closed vertices; however, they do not impact the difference $|V(\H)|-\pd(\H)$. 
 In fact, by Corollary \ref{BabyF}, $\pd(\H)=|V(\H)|-2$.

\begin{figure}[h] 
\caption{}\label{FL}
\begin{center}
\begin{tikzpicture} [thick, scale=0.6]
\shade [shading=ball, ball color=black] (4,-3) circle (.2);
\shade [shading=ball, ball color=black] (5,-2) circle (.2);
\path [pattern=north west lines, pattern  color=blue]   (4,-3)--(5,0)--(2.5,-3)--cycle;
\path [pattern=north west lines, pattern  color=blue]   (5,0)--(7,-1)--(5,-2)--cycle;

\path [pattern=north west lines, pattern  color=red]   (1,0)--(3,1)--(5,0)--cycle;

\path [pattern=north west lines, pattern  color=green]   (1.5,-1)--(2,-2)--(2,1)--cycle;

\draw  [shape=circle, color=red ] (1,0) circle (.2) node [left] {$\bf{f_1}$};
\draw  [shape=circle, color=green ] (1.5,-1) circle (.2) node [left] {$\bf{f_2}$};
\draw  [shape=circle, color=green] (2,-2) circle (.2) node [left] {$\bf{f_3}$};
\draw  [shape=circle, color=green] (2.5,-3) circle (.2) node [left] {$\bf{f_4}$};
\draw  [shape=circle, color=green] (3,-4) circle (.2) node [left] {$\bf{f_5}$};
\draw  [shape=circle, color=green ] (2,1.5) circle (.2) node [above] {$\bf{g_1}$};
\draw  [shape=circle, color=red ] (3,1) circle (.2) node [above] {$\bf{g_2}$};
\draw  [shape=circle, color=red] (4,0.5) circle (.2) node [above] {$\bf{g_3}$};
\draw  [shape=circle, color=red] (5,0) circle (.2) node [above] {$\bf{g_4}$};
\draw  [shape=circle, color=red ] (6,-0.5) circle (.2) node [above] {$\bf{g_5}$};
\draw  [shape=circle, color=red ] (7,-1) circle (.2) node [above] {$\bf{g_6}$};
\draw  [shape=circle, color=red ] (8,-1.5) circle (.2) node [above] {$\bf{g_7}$};

\draw [line width=1.2pt ] (1,0)--(2,1.5)  ;
\draw [line width=1.2pt , color=red   ] (1,0)--(3,1)  ;
\draw [line width=1.2pt, color=red ] (1,0)--(4,0.5)  ;
\draw [line width=1.2pt, color=red  ] (1,0)--(5,0)  ;
\draw [line width=1.2pt , color=red   ] (1,0)--(6,-0.5)  ;
\draw [line width=1.2pt , color=red   ] (1,0)--(7,-1)  ;
\draw [line width=1.2pt , color=red   ] (1,0)--(8,-1.5)  ;
\draw [line width=1.2pt, color=green  ] (1.5,-1)--(2,1.5)  ;
\draw [line width=1.2pt, color=green  ] (2,-2)--(2,1.5)  ;
\draw [line width=1.2pt, color=green  ] (2.5,-3)--(2,1.5)  ;
\draw [line width=1.2pt, color=green  ] (3,-4)--(2,1.5)  ;
\draw [line width=1.2pt ] (1.5,-1)--(3,1)  ;
\draw [line width=1.2pt ] (1.5,-1)--(4,0.5)  ;
\draw [line width=1.2pt ] (1.5,-1)--(5,0)  ;
\draw [line width=1.2pt ] (1.5,-1)--(6,-0.5)  ;
\draw [line width=1.2pt, color=green ] (2,-2)--(2,1)  ;
\draw [line width=1.2pt ] (2,-2)--(3,1)  ;
\draw [line width=1.2pt ] (2,-2)--(4,0.5)  ;
\draw [line width=1.2pt ] (2,-2)--(5,0)  ;
\draw [line width=1.2pt, color=green ] (2.5,-3)--(2,1)  ;
\draw [line width=1.2pt ] (2.5,-3)--(3,1)  ;
\draw [line width=1.2pt ] (2.5,-3)--(4,0.5)  ;
\draw [line width=1.2pt, color=green  ] (2.5,-3)--(3,-4)  ;
\draw [line width=1.2pt, color=green ] (3,-4)--(2,1)  ;
\draw [line width=1.2pt ] (3,-4)--(3,1)  ;

\draw [line width=1.2pt, color=red  ] (3,1)--(4,0.5)  ;

\end{tikzpicture}
\end{center}
\end{figure}

\end{Example}

\section{Projective Dimension of 1-Dimensional Hypergraphs}

A vertex $v$ in a 1-dimensional hypergraph $\H$ is called a {\em joint} if $\deg(v)\geq 3$. Let $v$ be a vertex in a hypergraph $\H$, and let $\H_1,\ldots,\H_r$ be the connected components of $\H_v$; if one of them, say $\H_1$, is a string hypergraph, we call $\H_1$ a {\it branch} of $\H$ (from $v$). 
This suggests the setting for the next result where we prove that if a hypergraph has a branch, then we can remove a few extremal vertices and keep track of the projective dimension.

\begin{Lemma} \label{ref2}
Let $\H$ be a hyerpgraph and $\BB$ be a branch of $\H$ it with at least $2$ vertices. Let $v_1$ be the endpoint of $\BB$, $v_2$ its neighbor, $v_3$ the neighbor of $v_2$( if there is one). Then 
\begin{itemize}
\item[$($a$)$] $\pd(\HH)=\pd(\HH_{v_1})+1$ 
 if $v_2$ is closed;
\item[$($b$)$] $\pd(\HH)=\pd(\HH_{v_1,v_2,v_3})+2$
 if $v_2$ is open.
\end{itemize}
\end{Lemma}

Assertion (b) generalizes \cite[Proposition~2.15]{LMa} to arbitrary hypergraphs.

\demo
Part (a) follows by Lemma \ref{red}. To prove (b) we apply Lemma \ref{ref} to $v_1$. Since $v_2$ is open, we have $\Q_{v_1}=\HH_{v_1,v_2,v_3}\cup{\{w\}}$ where $w$ is an isolated closed vertex (corresponding to the edge connecting $v_2$ and $v_3$), hence $\pd(\Q_{v_1})+1=\pd(\HH_{v_1,v_2,v_3})+2$. Now, Remark \ref{J+1} yields $\pd(\H_{v_1})\leq\pd(\HH_{v_1,v_2,v_3})+2$, therefore, by Lemma \ref{ref}, we have $\pd(\H)=\pd(\Q_{v_1})+1=\pd(\HH_{v_1,v_2,v_3})+2$.
\QED 
\bigskip

To study branches of hypergraphs, we need to recall that an {\it open string} is a string hypergraph $\H$ where every vertex is open except the two endpoints of $\H$ (which must be closed by separatedness). Every string consists of open strings which are (possibly) separated by closed vertices; more details on string hypergraphs and open strings can be found in \cite{LMa}.

Next, we define a more refined invariants of a string, which also keep track of the orientation. In fact, orientation appears to be crucial for branches inside general hypergraphs.
\begin{Definition}\label{W(S)}
Let $\mathcal{S}$ be a string and let $w=v_1$ and $v=v_n$ be its endpoints. Let $n_1,...,n_s$ be the number of opens in each open string in $\mathcal{S}$ starting from $w$, thus $n_s$ is the number of opens in the open string in $\SS$ attached to $v=v_n$.

(a) We say $\mathcal{S}$ is a {\rm 1-1 special configuration} if $\mathcal{S}$
does not contain two adjacent closed vertices and $n_{1}\equiv n_{s}\equiv 1$
mod $3$ and $n_{i}\equiv 2$ mod $3$ for $1<i<s$. The {\rm modularity from $v$} of a string $\mathcal{S}$ is the number $M(\mathcal{S};v)$ of pairwise disjoint 1-1 special configurations, counted starting from $v$. 

(b) We say $\mathcal{S}$ is a {\rm 1-0 special configuration from $v$} if $\mathcal{S}$
does not contain two adjacent closed vertices and $n_{s}\equiv 1$
mod $3$, $n_{1}\equiv 0$ mod 3 and $n_{i}\equiv 2$ mod $3$ for $1<i<s$.

(c) With $O(\mathcal{S};v)$ we denote the number of 1-0 special configurations with respect to $v$ which are disjoint both 1-1 special configurations and other 1-0 special configurations from $v$.

(d) We let $W(\mathcal{S};v)=|\{i|n_i\equiv 0\mbox{ mod }3\}|-O(\mathcal{S};v)$ be the number of open strings which have $3t$ open vertices for some $t\in \mathbb Z_+$ and are not part of a 1-0 special configuration with respect to $v$.

(e) The quotient of the division of $n-M(\SS;v)-W(\SS;v)$ by 3 is denoted by $q(\SS;v)$; the remainder of this division is denoted $\nr(\SS;v)$.
 \end{Definition}

In the next results we show that the number $\nr(\SS;v)$ essentially detects the point on a string (or branch) where we can cut the hypergraph without changing the projective dimension. In turn, this is the key point to find a simple way to compute the projective dimension of a number of 1-dimensional hypergraphs (see Theorem \ref{uniqueCycle}).

\begin{Remark}
The notions defined in Definition \ref{W(S)}.(b)-(e) are clearly sensitive to the choice of the orientation. For instance, if $v$ and $w$ are the two endpoints of $\SS$, one may have $W(\SS;v)\neq W(\SS;w)$. For example, if $\mathcal{S}$ is a string of length 7 with vertices $v_1,\ldots,v_7$ and $v_1,v_5,v_7$ are the only closed vertices of $\SS$, then $W(\mathcal{S};{v_1})=1$ whereas $W(\mathcal{S};{v_7})=0$.
\end{Remark}

If a 1-dimensional hypergraph $\H$ contains two adjacent closed vertices, then $\pd(\H)=\pd(\H')$ where $\H'$ is obtained by cancelling the edge connecting the two vertices (Proposition \ref{eq}). Thus, after replacing $\H$ by $\H'$,without loss of generality, we may assume in all the following statements that the 1-dimensional hypergraphs do not contain adjacent closed vertices.
 \bigskip

We can now prove the main technical result of this section. It gives a precise formula allowing us to detach all the strings from an arbitrary hypergraph. 
\begin{Proposition}\label{AttachString}
Let $\H$ be a hypergraph, $w$ a joint, and let $\SS$ be a branch of $\H$ from $w=v_0$ having vertices $v_1,...,v_n$ and containing no adjacent closed vertices.  Let $\ww{n}=\nr(\SS;v_n)$ and $q(\SS;v_n)$ be the numbers defined in Definition \ref{W(S)}.(e).  Let $E$ be the edge of $\SS$ 
\begin{itemize}
\item[$($i$)$] between the vertices  $v_{\widetilde{n}+1}$ and $v_{\widetilde{n}+2}$, if $v_1$ is open, the string of opens in $\SS$ ending in $w$ has $m\equiv 0$ (mod 3) open vertices, and $v_1$ is not part of 1-0 special configuration with respect to $v_n$; 

\item[$($ii$)$] between the vertices  $v_{\widetilde{n}}$ and  $v_{\widetilde{n}+1}$, in all other cases.
\end{itemize}

Then 
$$\pd(\H)=\pd(\H^E)=\pd(\widetilde{\H})+\pd(\mathcal{S'})$$
where $\H^E=\widetilde{\H}\cup \mathcal{S}'$ is obtained by cutting the edge $E$ of $\H$.  Moreover,  $$\pd(\mathcal{S}')=M(\mathcal{\SS};{v_n})+W(\mathcal{\SS};{v_n})+2q(\SS;v_n).$$

\end{Proposition}

\begin{Remark}
The proof of Proposition \ref{AttachString} has a subtle point, highlighted by the need to distinguish case (i) from all other cases. The following example illustrates it. Set $$\mathcal{H}=\{\{6\},\{6,0\},\{0,5\},\{5\},\{0,1\},\{1,2\},\{2,3\},\{3,4\},\{4\}\}$$ 
The hypergraph $\H$ contains the branch $\mathcal{S}=\{\{0,1\},\{1,2\},\{2,3\},\{3,4\},\{4\}\}$. One may use \cite{M2} to verify that $\pd(\H)=5$.

We have $W(\mathcal{S};v_{4})=1$ and $\widetilde{n}=4-1-3=0$. If we cut the edge between the vertices $v_{0}$ and $v_{1}$, we obtain two disjoint hypergraphs: a string with three closed vertices and an open string with four vertices. Then $\pd(\mathcal{H}^{E})=3+3=6$.
On the other hand, we have $\pd(\H)=5=\pd(\widetilde{\mathcal{H}})+\pd(\mathcal{S}')=3+2$
where $\widetilde{\mathcal{H}}$ and $\mathcal{S}'$ are obtained by cutting the edge between vertices $v_{1}$ and $v_{2}$ as in assumption (i) of Proposition \ref{AttachString}.
\end{Remark}

\noindent
{\bf Proof of Proposition \ref{AttachString}.} 
\noindent 
Let $\ww{n}=\nr(\SS;v_n)$ and $\ww{q}=q(\SS;v_n)$. We prove it by induction on the number of open strings of $\mathcal{S}$. We may assume $n>2$. The inductive step is the same for both (i) and (ii), thus we first prove the inductive step and the base case later.

In the induction step, we may assume $\mathcal{S}$ has at least two open strings. Let $\SS_s$ be the open string in
$\SS$, having $v_n$ as an endpoint; let us denote by $v_{s1}$ the other endpoint of $\SS_s$ and by $q_s$ the quotient obtained when we divide the number $n_s$ of opens in $\SS_s$ by 3. If $n_s=3q_s$ then $q_s$ iterations of Lemma \ref{ref2}.(b) and one iteration of Lemma \ref{ref2}.(a) yield $\pd(\mathcal{H})=\pd(\mathcal{H}')+1+2q_s$, where $\mathcal{H}'$ is the hypergraph obtained by removing from $\H$ the last $3q_s+1$ vertices of $\SS$. We denote by $\SS''$ the string in $\H'$ obtained after this procedure has removed $\SS_s$ from $\SS$; also, let $v_1,\ldots,v_{n''}$ be the vertices in $\SS''$ let $q''=q(\SS'';v_{n''})$ (note that $n''=n-3q_s-1$ and $q''=\ww{q}-q_s$) and $\ww{n''}=\nr(\SS'';v_{n''}$). Clearly $\SS''$ has one less open string than $\mathcal{S}$; moreover, $W(\mathcal{S''};{v_{s1}})=W(\mathcal{S};{v_n})-1$ and $M(\mathcal{S''};{v_{s1}})=M(\mathcal{S};{v_n})$. 
Then, applying the induction hypothesis to $\H'$ (and the string $\SS''$ in it), we obtain $$\pd(\mathcal{H}')=\pd(\ww{\H'})+\pd(\SS''')\;\;\; \mbox{ and }\;\;\;\pd(\mathcal{S}''')=M(\mathcal{S''};{v_{s1}})+W(\mathcal{S''};{v_{s1}})+2q'',$$
where $\ww{\H'}$ and $\SS'''$ are the disjoint components of $(\H')^{E'}$.
Notice that $\ww{n}=n-M(\mathcal{S};{v_n})-W(\mathcal{S};{v_n})-3\ww{q}=n''-M(\mathcal{S}'';{v_{s1}})-W(\mathcal{S}'';{v_{s1}})-
3q''=\widetilde{n}''$. Since $\ww{n}=\ww{n}''$, then $E=E'$ and $\ww{\H}=\ww{\H'}$. Since $\pd(\H)=2q_s+1+\pd(\H')$, we obtain
$$\pd(\H)=2q_s+1+\pd(\ww{\H})+M(\mathcal{S''};{v_{s1}})+W(\mathcal{S''};{v_{s1}})+2q''$$
Therefore $\pd(\mathcal{H})=\pd(\ww{\H})+M(\mathcal{S};{v_n})+W(\mathcal{S};{v_n})+2\ww{q}$. The case where $\mathcal{S}_{s}$ has $n_s=3q_s+2$ open vertices follows similarly, after $q_s+1$ applications of Lemma \ref{ref2}. 

If $n_s=3q_s+1$ we have three scenarios to consider: when $\mathcal{S}_{s}$ is part of a 1-1 special configuration; when $\mathcal{S}_{s}$ is part of a 1-0 special configuration;  when $\mathcal{S}_{s}$ is neither part of a 1-1 special configuration nor of a 1-0 special configuration. In this latter case, $\mathcal{S}$ is a string with $n_i\equiv 2$ mod $3$ for $i=1,...,s-1$, and $n_s\equiv 1$, and $0=W(\mathcal{S};{v_n})=M(\mathcal{S};{v_n})$. When $v_1$ is closed, we have $\widetilde{n}=0$, and when $v_1$ is part of the open vertex of the $n_1$, we have $\widetilde{n}=2$.  Repeated applications of Lemma \ref{ref2} yield $\pd(\mathcal{S'})=2\left\lfloor \frac{n}{3}\right\rfloor=2\ww{q}$.

When $\mathcal{S}_s$ is part of a 1-1 special configuration $\mathcal{M}\subseteq \SS$, we remark that, by definition, $\mathcal{M}$ has $3l+2$ vertices (for some integer $l$) and endpoints $v_n$ and, say, $v'$. Applying Lemma \ref{ref2}.(b) for $l$ times, and Lemma \ref{ref2}.(a) once, one obtains $\pd(\mathcal{H})=\pd(\mathcal{H}')+2l+1$. Then we can apply induction on $\mathcal{H}'$, and a proof similar to the above (with $W(\mathcal{S}';{v'})=W(\mathcal{S};{v_n})$ and $M(\mathcal{S}';{v'})=M(\mathcal{S};{v_n})-1$) gives the conclusion.

When  $\mathcal{S}_s$ is part of 1-0 special configuration $\mathcal{O}$, we observe that by definition the string $\mathcal{O}$ has $3l+1$ vertices (for some $l$),  and endpoints $v_n$ and, say, $v'$. We can then apply Lemma \ref{ref2}.(b) for $l$ times to obtain $\pd(\mathcal{H})=\pd(\mathcal{H}')+2l$. Then, the induction hypothesis applied to $\mathcal{H}'$ and a proof similar to the above (with $W(\mathcal{S}';{v'})=W(\mathcal{S};{v_n})$ and $M(\mathcal{S}';{v'})=M(\mathcal{S};{v_n})$) yield the final statement.

Observe that the inductive step does not change the structure of the first open string of $\SS$. We may now prove the base case, where $\SS$ is a string of opens. There are three cases to consider: $n=3q$, $n=3q+1$, and $n=3q+2$ for some integer $q$. In the first case we apply Lemma \ref{ref2}.(b) $q$ times to obtain the conclusion. Notice that in this case, $W(\mathcal{S};{v_n})=0$, $M(\mathcal{S};{v_n})=0$, $\ww{n}=0$, and $\ww{q}=q(\SS;v_n)=q$.

Assume $n=3q+1$. If $v_1$ is open then we are under assumption (i), and $W(\mathcal{S};{v_n})=1$ and $\ww{n}=0$. We apply Lemma \ref{ref2}.(b) $q$ times and obtain 
$$\mbox{pd}(\mathcal{H})=\mbox{pd}\widetilde{\mathcal{H}}+2q,$$
where $\widetilde{\mathcal{H}}$ is as in the assumptions, because we cut the edge between $v_{\ww{n}+1}$ and $v_{{\ww{n}}+2}$ when we apply Lemma \ref{ref2}.(b) $q$ times. 
On the other hand, if $v_1$ is closed, then  $W(\mathcal{S};{v_n})=0$ and $\ww{n}=1$, because $\SS$ has $3q-1$ open vertices and $v_1$ and $v_n$ are closed. This time $q$ applications of Lemma \ref{ref2}.(b) give the conclusion. 

We may then assume $n=3q+2$. If $v_1$ is open, then $W(\mathcal{S};{v_n})=0$ and $\ww{n}=2$. By applying Lemma \ref{ref2}.(b) $q$ times, we obtain  
$$\mbox{pd}(\mathcal{H})=\mbox{pd}\widetilde{\mathcal{H}}+2q,$$
as before. Instead, when $v_1$ is closed, then $\SS$ has $3q$ open vertices and  $W(\mathcal{S};{v_n})=1$ and $\ww{n}=1$. After applying Lemma \ref{ref2}.(b) $q$ times, we apply Lemma \ref{ref2}.(a) once. We obtain   $$\mbox{pd}(\mathcal{H})=\mbox{pd}\widetilde{\mathcal{H}}+2q+1,$$
and we cut the edge between $v_1$ and $v_2$ when we apply  Lemma \ref{ref2}.(a). This finishes the proof.
\QED 
\bigskip

Notice that the branch $\ww{S}$ of $S$ in $\ww{\H}$ has at most two vertices. Moreover, if it has two vertices, the vertex connected to the joint $w$ must be open.  The next step consists in finding the projective dimension of $\ww{\H}$ when each branch of the hypergraph has length at most $2$. The following proposition gives a reduction that detach the branches of the hypergraph in a controlled way.

\begin{Proposition}\label{cutTree}
Let $\H$ be a 1-dimensional hypergraph, $w$ a joint in $\H$, $\SS$ a branch departing from $w$  with $n$ vertices and containing no adjacent closed vertices, and let $E$ be the edge connecting $w$ to $\SS$. Then $$\pd(\H)=\pd(\H')$$
where $\H'$ is the following hypergraph:
\begin{itemize}

\item[$($a$)$] if $n=1$, then $\H'=\H:E$, i.e. $\H'$ is obtained by  cancelling $E$;

\item[$($b$)$] if $n=2$, 
then $\H'=\H_w$, i.e. $\H'$ is obtained by removing $w$.
\end{itemize}
\end{Proposition}

\demo
For assertion (a), let $v_1$ be the vertex in $\SS$, we apply Lemma \ref{ref} to obtain 
$$\pd(\H)=\max\{\pd(\H_{v_1}),\pd(\Q_{v_1})+1\}$$
Observe that $\Q_{v_1}=(\H:E)\,\backslash {\{v_1\}}=\H'\backslash {\{v_1\}}$ and $\H_{v_1}=\H'\backslash {\{v_1\}}\cup \{w\}$ (if $\{w\}\in \H$, i.e. if $w$ is already closed in $\H$, then $\H_{v_1}=\H'\backslash {\{v_1\}}$). Then, by Lemma \ref{H^0}, $\pd(\H_{v_1})\leq \pd(\Q_{v_1})+1$, thus $\pd(\H)=\pd(\Q_{v_1})+1=\pd(\H')$.

For part (b), let $v_1$ and $v_2$ be the two vertices of $\SS$ and assume $v_1$ is the neighbor of $w$. By Lemma \ref{ref} we have 
$$\pd(\H)=\max\{\pd(\H_{v_2}),\pd(\Q_{v_2})+1\}$$
Notice that $\SS$ becomes a branch of length 1 in $\H_{v_2}$ and then, by (a), $\pd(\H_{v_2})=\pd(\H_{v_2}:E)$. Since $v_1$ is open (as remarked before Proposition \ref{cutTree}), then $\H_{v_2}:E=\H:E$; on the other hand, one has $\Q_{v_2}=\H_{v_2}:v_2=\H_w\backslash {\{v_2\}}$. Thus $\pd(\Q_{v_2})+1=\pd(\H_w)=\pd(\H')$ and it suffices to show that $\pd(\H:E)\leq \pd(\H')$. Finally, since the branch $\SS$ in $\H'$ is a string of length $2$, then  $\pd((\H:E)_w)+1=\pd(\H')$. Now, Remark \ref{J+1} gives $\pd(\H:E)\leq \pd((\H:E)_w)+1$, from which the conclusion follows.
\QED 
\bigskip

We can now prove the main result of this section, stating that we have a simple procedure to compute the projective dimension of many 1-dimensional hypergraphs.
\begin{Theorem}\label{uniqueCycle}
Let $\H$ be a 1-dimensional hypergraph. If each of its connected components contains at most one cycle, then by using the reduction process  of Proposition \ref{AttachString} and \ref{cutTree}, one can obtain $\pd(\H)$.
\end{Theorem}

\demo
It is enough to prove the statement when $\H$ is a connected hypergraph having at most one cycle. We may further assume that no closed vertices of $\H$ are adjacent, because otherwise we can cancel the edge connecting them without changing the projective dimension -- by Proposition \ref{eq}. Now, if $\H$ is a cycle, then the statement follows by \cite[Theorem~3.4]{LMa}. If $\H$ contains no cycles, then iterated applications of Propositions  \ref{AttachString} and \ref{cutTree} allow us to replace $\H$ by a disjoint union of strings. If $\H$ strictly contains a cycle, then by assumption $\H$ is a cycle with one or more trees attached to its vertices, and by repeated applications of Propositions  \ref{AttachString} and \ref{cutTree} we may replace $\H$ by a disjoint union of strings and, possibly, one cycle. 
In each of these scenarios, the main theorems in \cite{LMa} now allow us to compute the projective dimension of each component, and, therefore, the projective dimension of the original hypergraph $\H$.
\QED 
\bigskip

In the Appendix we implemented explicitly two algorithmic procedures that can be  employed to compute $\pd(\H)$ (in particular, see Algorithm \ref{AlgPd}).

\begin{Remark}\label{morecycles}
$($a$)$ Theorem \ref{uniqueCycle} may also be applied in certain cases where the connected components of $\H$ contain more than one cycle, provided that all cycles except, possibly, one per connected component, at the end of the reduction process of Proposition \ref{cutTree} become either strings, or are pairwise disjoint, or a combination of these two possibilities. This situation appears fairly frequently,  because if a joint has a branch $\SS$ with $\ww{n}=2$ (in the statement of Proposition \ref{AttachString}), then by Propositions \ref{AttachString} and \ref{cutTree}.(b) we can remove the joint; the removal unfolds the cycle and makes it become a string whenever the joint is one of the vertices of the cycle. \\
$($b$)$ Any 1-dimensional hypergraph is obtained by attaching together any number of trees and cycles. Since Theorem \ref{uniqueCycle} applies to all 1-dimensional hypergraphs $\H$ having an arbitrary number of trees and at most one cycle and to some cases where $\H$ has multiple cycles cases (see (a)), then Theorem \ref{uniqueCycle} provides an effective method to determine the projective dimension of a wide class of 1-dimensional hypergraphs.
\end{Remark}
The next example illustrates Remark \ref{morecycles} in a concrete situation where $\H$ contains 3 different cycles.

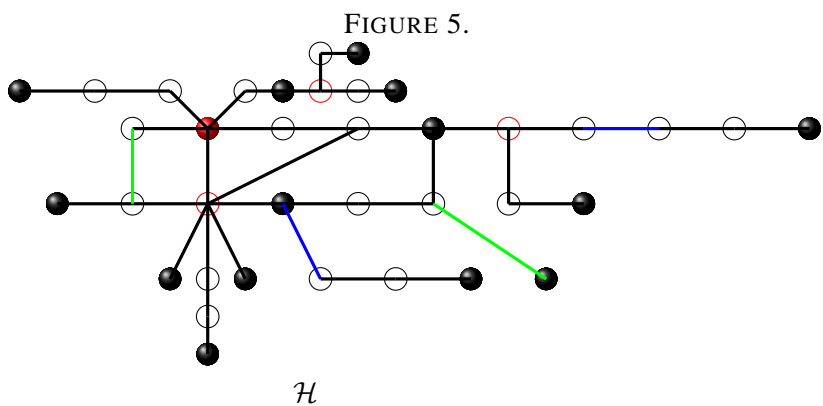
\begin{figure}[h] 
\caption{}\label{cutH}
\begin{center}
\begin{tikzpicture}

\shade [shading=ball, ball color=black]  (0,0) circle (.15);
\draw  [shape=circle] (1,0) circle (.15);
\draw  [shape=circle, color=red] (2,0) circle (.15);
\shade [shading=ball, ball color=black]  (3,0) circle (.15);
\draw  [shape=circle] (4,0) circle (.15);
\draw  [shape=circle] (5,0) circle (.15);
\draw  [shape=circle] (6,0) circle (.15);
\shade [shading=ball, ball color=black]  (7,0) circle (.15);

\shade [shading=ball, ball color=black] (1.5,-1) circle (.15);
\draw  [shape=circle] (2,-1) circle (.15);
\draw  [shape=circle] (2,-1.5) circle (.15);
\shade [shading=ball, ball color=black] (2,-2) circle (.15);
\shade [shading=ball, ball color=black]  (2.5,-1) circle (.15);
\draw  [shape=circle] (3.5,-1) circle (.15);
\draw  [shape=circle] (4.5,-1) circle (.15);
\shade [shading=ball, ball color=black]  (5.5,-1) circle (.15);
\shade [shading=ball, ball color=black]  (6.5,-1) circle (.15);
\draw  [shape=circle] (1,1) circle (.15);
\shade [shading=ball, ball color=red] (2,1) circle (.15);
\draw  [shape=circle] (3,1) circle (.15);
\draw  [shape=circle] (4,1) circle (.15);
\shade [shading=ball, ball color=black]  (5,1) circle (.15);
\draw  [shape=circle, color=red] (6,1) circle (.15);
\draw  [shape=circle]  (7,1) circle (.15);
\draw  [shape=circle]  (8,1) circle (.15);
\draw  [shape=circle]  (9,1) circle (.15);
\shade [shading=ball, ball color=black]   (10,1) circle (.15);

\draw  [shape=circle] (1.5,1.5) circle (.15);
\draw  [shape=circle] (0.5,1.5) circle (.15);
\shade [shading=ball, ball color=black]  (-0.5,1.5) circle (.15);
\draw  [shape=circle] (2.5,1.5) circle (.15);
\shade [shading=ball, ball color=black]  (3,1.5) circle (.15);
\draw  [shape=circle, color=red] (3.5,1.5) circle (.15);
\draw  [shape=circle] (4,1.5) circle (.15);
\shade [shading=ball, ball color=black] (4.5,1.5) circle (.15);
\draw  [shape=circle] (3.5,2) circle (.15);
\shade [shading=ball, ball color=black]  (4,2) circle (.15);

\draw [line width=1.2pt ] (0,0)--(1,0)  ;
\draw [line width=1.2pt ] (1,0)--(2,0)  ;
\draw [line width=1.2pt ] (2,0)--(3,0)  ;
\draw [line width=1.2pt ] (4,0)--(3,0)  ;
\draw [line width=1.2pt ] (4,0)--(5,0)  ;
\draw [line width=1.2pt ] (6,0)--(7,0)  ;
\draw [line width=1.2pt ] (6,0)--(6,1)  ;
\draw [line width=1.2pt ] (9,1)--(10,1)  ;
\draw [line width=1.2pt ] (8,1)--(9,1)  ;
\draw [line width=1.2pt, color=blue ] (7,1)--(8,1)  ;
\draw [line width=1.2pt ] (7,1)--(6,1)  ;
\draw [line width=1.2pt ] (5,1)--(6,1)  ;
\draw [line width=1.2pt ] (5,1)--(5,0)  ;
\draw [line width=1.2pt ] (5,1)--(4,1)  ;
\draw [line width=1.2pt ] (3,1)--(4,1)  ;
\draw [line width=1.2pt ] (2,1)--(3,1)  ;
\draw [line width=1.2pt ] (1,1)--(2,1)  ;
\draw [line width=1.2pt ] (2,1)--(2,0)  ;
\draw [line width=1.2pt ] (2,0)--(4,1)  ;
\draw [line width=1.2pt, color=green ] (1,1)--(1,0)  ;
\draw [line width=1.2pt ] (2,1)--(1.5,1.5)  ;
\draw [line width=1.2pt ] (0.5,1.5)--(1.5,1.5)  ;
\draw [line width=1.2pt ] (-0.5,1.5)--(0.5,1.5)  ;
\draw [line width=1.2pt ] (2,1)--(2.5,1.5)  ;
\draw [line width=1.2pt ] (3,1.5)--(2.5,1.5)  ;
\draw [line width=1.2pt ] (3,1.5)--(3.5,1.5)  ;
\draw [line width=1.2pt ] (4,1.5)--(3.5,1.5)  ;
\draw [line width=1.2pt ] (4,1.5)--(4.5,1.5)  ;
\draw [line width=1.2pt ] (3.5,1.5)--(3.5,2)  ;
\draw [line width=1.2pt ] (4,2)--(3.5,2)  ;
\draw [line width=1.2pt ] (2,-1)--(2,0)  ;
\draw [line width=1.2pt ] (2.5,-1)--(2,0)  ;
\draw [line width=1.2pt ] (1.5,-1)--(2,0)  ;
\draw [line width=1.2pt ] (2,-1)--(2,-1.5)  ;
\draw [line width=1.2pt ] (2,-1.5)--(2,-2)  ;
\draw [line width=1.2pt, color=blue] (3.5,-1)--(3,0)  ;
\draw [line width=1.2pt ] (3.5,-1)--(4.5,-1)  ;
\draw [line width=1.2pt ] (4.5,-1)--(5.5,-1)  ;
\draw [line width=1.2pt, color=green ] (6.5,-1)--(5,0)  ;
\path (3,-2.5)--(3,-2.5) node [pos=.5, right ] {$\H$};

\end{tikzpicture}
\end{center}
\end{figure}

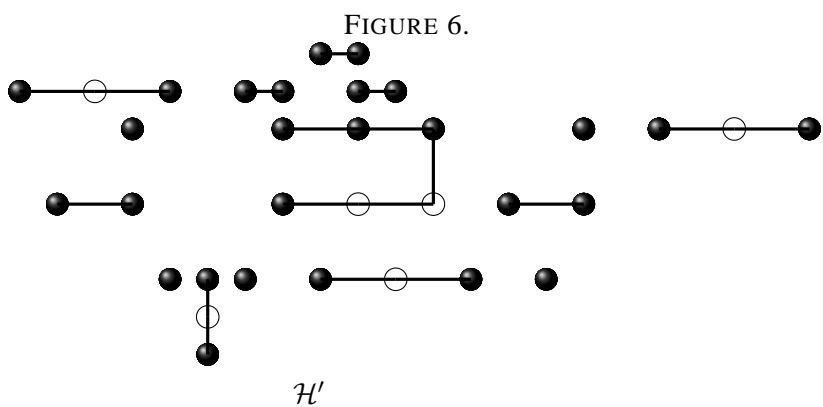
\begin{figure}[h] 
\caption{}\label{cutH'}
\begin{center}
\begin{tikzpicture}

\shade [shading=ball, ball color=black]  (0,0) circle (.15);
\shade [shading=ball, ball color=black]  (1,0) circle (.15);

\shade [shading=ball, ball color=black]  (3,0) circle (.15);
\draw  [shape=circle] (4,0) circle (.15);
\draw  [shape=circle] (5,0) circle (.15);
\shade [shading=ball, ball color=black]  (6,0) circle (.15);
\shade [shading=ball, ball color=black]  (7,0) circle (.15);

\shade [shading=ball, ball color=black] (1.5,-1) circle (.15);
\shade [shading=ball, ball color=black]  (2,-1) circle (.15);
\draw  [shape=circle] (2,-1.5) circle (.15);
\shade [shading=ball, ball color=black] (2,-2) circle (.15);
\shade [shading=ball, ball color=black]  (2.5,-1) circle (.15);
\shade [shading=ball, ball color=black]  (3.5,-1) circle (.15);
\draw  [shape=circle] (4.5,-1) circle (.15);
\shade [shading=ball, ball color=black]  (5.5,-1) circle (.15);
\shade [shading=ball, ball color=black]  (6.5,-1) circle (.15);
\shade [shading=ball, ball color=black]  (1,1) circle (.15);

\shade [shading=ball, ball color=black] (3,1) circle (.15);
\shade [shading=ball, ball color=black]  (4,1) circle (.15);
\shade [shading=ball, ball color=black]  (5,1) circle (.15);

\shade [shading=ball, ball color=black]   (7,1) circle (.15);
\shade [shading=ball, ball color=black]   (8,1) circle (.15);
\draw  [shape=circle]  (9,1) circle (.15);
\shade [shading=ball, ball color=black]   (10,1) circle (.15);

\shade [shading=ball, ball color=black]  (1.5,1.5) circle (.15);
\draw  [shape=circle] (0.5,1.5) circle (.15);
\shade [shading=ball, ball color=black]  (-0.5,1.5) circle (.15);
\shade [shading=ball, ball color=black]  (2.5,1.5) circle (.15);
\shade [shading=ball, ball color=black]  (3,1.5) circle (.15);

\shade [shading=ball, ball color=black]  (4,1.5) circle (.15);
\shade [shading=ball, ball color=black] (4.5,1.5) circle (.15);
\shade [shading=ball, ball color=black]  (3.5,2) circle (.15);
\shade [shading=ball, ball color=black]  (4,2) circle (.15);

\draw [line width=1.2pt ] (0,0)--(1,0)  ;

\draw [line width=1.2pt ] (4,0)--(3,0)  ;
\draw [line width=1.2pt ] (4,0)--(5,0)  ;
\draw [line width=1.2pt ] (6,0)--(7,0)  ;

\draw [line width=1.2pt ] (9,1)--(10,1)  ;
\draw [line width=1.2pt ] (8,1)--(9,1)  ;

\draw [line width=1.2pt ] (5,1)--(5,0)  ;
\draw [line width=1.2pt ] (5,1)--(4,1)  ;
\draw [line width=1.2pt ] (3,1)--(4,1)  ;

\draw [line width=1.2pt ] (0.5,1.5)--(1.5,1.5)  ;
\draw [line width=1.2pt ] (-0.5,1.5)--(0.5,1.5)  ;

\draw [line width=1.2pt ] (3,1.5)--(2.5,1.5)  ;

\draw [line width=1.2pt ] (4,1.5)--(4.5,1.5)  ;

\draw [line width=1.2pt ] (4,2)--(3.5,2)  ;

\draw [line width=1.2pt ] (2,-1)--(2,-1.5)  ;
\draw [line width=1.2pt ] (2,-1.5)--(2,-2)  ;

\draw [line width=1.2pt ] (3.5,-1)--(4.5,-1)  ;
\draw [line width=1.2pt ] (4.5,-1)--(5.5,-1)  ;

\path (3,-2.5)--(3,-2.5) node [pos=.5, right ] {$\H'$};

\end{tikzpicture}
\end{center}
\end{figure}

\begin{Example}
The hypergraph $\H'$ in Figure \ref{cutH'} is obtained by applying Propositions \ref{AttachString} and \ref{cutTree} repeatedly from the hypergraph $\H$ in Figure \ref{cutH}: we cut the blue faces, cancel the green faces and remove the red vertices. Our procedure gives $\pd(\H)=\pd(\H')=28$.
\end{Example}

We now define stars and use them to introduce a more complicated class of hypergraph, obtained by connecting together stars via their centers. We will provide explicit combinatorial formulas for their projective dimensions.
\begin{Definition}\label{openstar}
A connected hypergraph $\H$ is called a {\em star} if either $|V(\H)|=1$ (and we call the only vertex of $\H$ its {\em  center}) or $\H$ does not contain any cycle, no adjacent closed vertices, and it contains precisely one joint (called the {\em center} of the star). $\H$ is called a $d$-star if $\H$ is a star and every branch in $\H$ has length at most $d$.

An {\em open star} ({\em closed star}, resp.) is a star whose center is an open (closed, resp.) vertex.
\end{Definition}

Note that $\H$ is a $0$-star if and only if $|V(\H)|=1$. Also, any $(d-1)$-star is also a $d$-star when $d>1$; thus we say that $\H$ is a {\em proper $d$-star} if $\H$ is a $d$-star and $\H$ is not a $(d-1)$-star (i.e. if $\H$ contains at least one branch of length $d$). We now give a few more definitions, which can be interpreted as natural generalizations of strings and cycle hypergraphs to stars. The only exception is that for string hypergraphs the assumption of separatedness forces the endpoints of the string to be closed vertices, whereas for strings of stars this need not be the case:

\begin{Definition}\label{stringcycle}
 A {\rm string (cycle, tree, resp.) of stars} is a hypergraph $\H$ consisting of a (finite) collection of stars where each star is only connected to other stars via its center and the centers of the stars form a string (cycle, tree, resp.). \end{Definition}

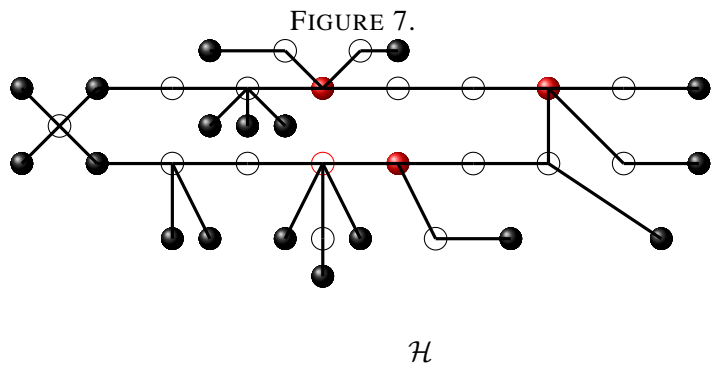
\begin{figure}[h] 
\caption{}\label{Estar}
\begin{center}
\begin{tikzpicture}
\shade [shading=ball, ball color=black] (-2,1) circle (.15);
\shade [shading=ball, ball color=black] (-2,0) circle (.15);
\shade [shading=ball, ball color=black] (0.5,-1) circle (.15);
\shade [shading=ball, ball color=black] (0,-1) circle (.15);
\shade [shading=ball, ball color=black] (-1,1) circle (.15);
\shade [shading=ball, ball color=black] (-1,0) circle (.15);
\shade [shading=ball, ball color=black] (0.5,0.5) circle (.15);
\shade [shading=ball, ball color=black] (1.5,0.5) circle (.15);
\shade [shading=ball, ball color=black] (1,0.5) circle (.15);
\draw  [shape=circle] (-1.5,0.5) circle (.15);
\draw  [shape=circle] (0,0) circle (.15);
\draw  [shape=circle] (0,1) circle (.15);
\draw  [shape=circle] (1,0) circle (.15);
\draw  [shape=circle, color=red] (2,0) circle (.15);
\shade [shading=ball, ball color=red]  (3,0) circle (.15);
\draw  [shape=circle] (4,0) circle (.15);
\draw  [shape=circle] (5,0) circle (.15);
\draw  [shape=circle] (6,0) circle (.15);
\shade [shading=ball, ball color=black]  (7,0) circle (.15);

\shade [shading=ball, ball color=black] (1.5,-1) circle (.15);
\draw  [shape=circle] (2,-1) circle (.15);
\shade [shading=ball, ball color=black] (2,-1.5) circle (.15);

\shade [shading=ball, ball color=black]  (2.5,-1) circle (.15);
\draw  [shape=circle] (3.5,-1) circle (.15);
\shade [shading=ball, ball color=black]  (4.5,-1) circle (.15);

\draw  [shape=circle] (1,1) circle (.15);
\shade [shading=ball, ball color=red] (2,1) circle (.15);
\draw  [shape=circle] (3,1) circle (.15);
\draw  [shape=circle] (4,1) circle (.15);
\shade [shading=ball, ball color=red]  (5,1) circle (.15);
\draw  [shape=circle] (6,1) circle (.15);
\shade [shading=ball, ball color=black]  (7,1) circle (.15);
\draw  [shape=circle] (1.5,1.5) circle (.15);
\shade [shading=ball, ball color=black] (0.5,1.5) circle (.15);

\draw  [shape=circle] (2.5,1.5) circle (.15);
\shade [shading=ball, ball color=black]  (3,1.5) circle (.15);
\shade [shading=ball, ball color=black]  (6.5,-1) circle (.15);
\draw [line width=1.2pt ] (-1.5,0.5)--(-2,0)  ;
\draw [line width=1.2pt ] (-1.5,0.5)--(-2,1)  ;
\draw [line width=1.2pt ] (0,0)--(0.5,-1)  ;
\draw [line width=1.2pt ] (0,0)--(0,-1)  ;
\draw [line width=1.2pt ] (-1.5,0.5)--(-1,0)  ;
\draw [line width=1.2pt ] (-1.5,0.5)--(-1,1)  ;
\draw [line width=1.2pt ] (0,0)--(-1,0)  ;
\draw [line width=1.2pt ] (0,1)--(-1,1)  ;
\draw [line width=1.2pt ] (0,1)--(1,1)  ;
\draw [line width=1.2pt ] (0,0)--(1,0)  ;
\draw [line width=1.2pt ] (1,0)--(2,0)  ;
\draw [line width=1.2pt ] (2,0)--(3,0)  ;
\draw [line width=1.2pt ] (4,0)--(3,0)  ;
\draw [line width=1.2pt ] (4,0)--(5,0)  ;
\draw [line width=1.2pt ] (6,0)--(7,0)  ;
\draw [line width=1.2pt ] (6,0)--(5,1)  ;

\draw [line width=1.2pt ] (7,1)--(6,1)  ;
\draw [line width=1.2pt ] (5,1)--(6,1)  ;
\draw [line width=1.2pt ] (1,1)--(.5,0.5)  ;
\draw [line width=1.2pt ] (1,1)--(1.5,0.5)  ;
\draw [line width=1.2pt ] (1,1)--(1,0.5)  ;
\draw [line width=1.2pt ] (5,1)--(5,0)  ;
\draw [line width=1.2pt ] (5,1)--(4,1)  ;

\draw [line width=1.2pt ] (3,1)--(4,1)  ;
\draw [line width=1.2pt ] (2,1)--(3,1)  ;
\draw [line width=1.2pt ] (1,1)--(2,1)  ;

\draw [line width=1.2pt ] (2,1)--(1.5,1.5)  ;
\draw [line width=1.2pt ] (0.5,1.5)--(1.5,1.5)  ;

\draw [line width=1.2pt ] (2,1)--(2.5,1.5)  ;
\draw [line width=1.2pt ] (3,1.5)--(2.5,1.5)  ;
\draw [line width=1.2pt ] (2,-1)--(2,0)  ;
\draw [line width=1.2pt ] (2.5,-1)--(2,0)  ;
\draw [line width=1.2pt ] (1.5,-1)--(2,0)  ;
\draw [line width=1.2pt ] (2,-1)--(2,-1.5)  ;

\draw [line width=1.2pt] (3.5,-1)--(3,0)  ;
\draw [line width=1.2pt ] (3.5,-1)--(4.5,-1)  ;

\draw [line width=1.2pt] (6.5,-1)--(5,0)  ;
\path (3,-2.5)--(3,-2.5) node [pos=.5, right ] {$\H$};

\end{tikzpicture}
\end{center}
\end{figure}

Strings, cycles or trees of $d$-stars can also be thought as being obtained by taking a string (cycle, tree, resp.) hypergraph and attaching to some (or all) of its vertices strings of length at most $d$.
\begin{Example}
The hypergraph depicted in Figure \ref{Estar} illustrates an example of a cycle of $2$-stars.
\end{Example}

First, we prove a simple formula for the projective dimension of open strings of ``small'' stars (i.e. $d$-stars with $d\leq 2$).
\begin{Proposition}
If $\H$ is a disjoint union of trees and cycles of proper $2$-stars, then $\pd(\H)=|V(\H)|-T(\H)$ where $T(\H)$ is the number of proper $2$-stars in $\H$.
\end{Proposition}

\demo
By Proposition \ref{eq} we may assume there are no adjacent closed vertices. Also, by assumption, every star has at least a branch of length two, thus we can apply Proposition \ref{cutTree}.(b) to each star and obtain 
$$\pd(\H)=\pd(\H')$$
where $\H'$ is obtained from $\H$ by removing the center of each star. Since every star in $\H$ is a two star, then any open vertex in $\H$ is adjacent to the center of a 2-star, thus it becomes closed after removing the center. Therefore, all vertices in $\H'$ are now closed, i.e. $\H'$ is a saturated hypergraph with $|V(\H)|-T(\H)$ vertices. Therefore, $\pd(\H)=|V(\H)|-T(\H)$.
\QED
 \bigskip
 
In the study of string and cycle hypergraphs $\H$, it was introduced a purely combinatorial invariant, called {\it modularity} \cite[Definition~3.1]{LMa}, employed in the formulas for $\pd(\H)$. Here we naturally generalize this concept to strings and cycles of stars.
\begin{Definition}
We say that a string of stars $\H$ is a {\em 1-1 special star configuration} if $\H$ does not contain
 two adjacent closed stars and $n_{1}\equiv n_{s}\equiv 1$ mod $3$
 and $n_{i}\equiv 2$ mod $3$ for $1<i<s$, where $n_i$ is the number of open stars in the $i$-th open string of star in $\H$. The {\em star modularity} $M^*(\H)$ of a string or cycle of stars of $\H$ is the maximal number of pairwise of
 disjoint 1-1 special star configurations contained in $\H$. Similarly, the star modularity of a disjoint union of strings or cycles of stars is the sum of the star modularity of each connected component.
\end{Definition}

The following proposition has a similar flavour as \cite[Theorems~3.4, 4.3]{LMa} and it provides an effective combinatorial formula to compute the projective dimension of strings and cycles of 2-stars.

\begin{Proposition}\label{EstarsAndStars}
Let $\H$ be a string or cycle of $2$-stars. Let $T(\H)$ be the number of proper $2$-stars, let $\H^*$ be obtained from $\H$ by removing the centers of the proper 2-stars, and let $s^*(\H^{*})$ be the number of open strings of stars in $\H^{*}$. 
\begin{enumerate}
\item If $\H$ is an open cycle with $n_1$ 1-stars, then $\pd(\H)=|V(\H)|-1-\left\lfloor \frac{n_1-2}{3}\right\rfloor $;
\item in all other cases $\pd(\H)=|V(\H)|-T(\H)-s^*(\H^{*})-\sum_{i=1}^{s}\left\lfloor \frac{n_{i}-1}{3}\right\rfloor +M^*(\H^{*})$, where $n_{i}$ is the number of open stars in each open string of stars in $\H^{*}$. 
\end{enumerate} 
\end{Proposition}

\demo
Applying Proposition \ref{cutTree}.(b) to all proper 2-stars we obtain $\pd(\H)=\pd(\H^*)$, and observe that $\H^*$ is the disjoint union of strings or cycles of 1-stars. Let $\H^{**}$ be the hypergraph obtained after applying Proposition \ref{cutTree}.(a) to each proper 1-star. By definition of cancellation (Definition \ref{deform}.(iii)), $\H^{**}$ has the same modularity and number of open strings as $\H^*$, and each open string has the same number of vertices as $\H^*$, because $\H^{**}$ is obtained by cancelling the edges of the branches of 1-stars. 
Then $\H^{**}$ is a disjoint union of closed vertices and string or cycle hypergraphs, with $s^*(\H^*)$ open strings and $M^*(\H^{**})=M^*(\H^*)$. Now the conclusion follows by applying \cite[Theorems~3.4,~4.3]{LMa} to each connected component.
\QED 
\bigskip

\begin{Example}
Figure \ref{Estar} below depicts a cycle of 2-stars $\H$; Figure \ref{EstarReduce} shows the hypergraph $\H^*$ as defined in Proposition \ref{EstarsAndStars} (obtained by removing the red-labelled vertices). Then $\pd(\H)=37-4-3+1=31$.
\end{Example}

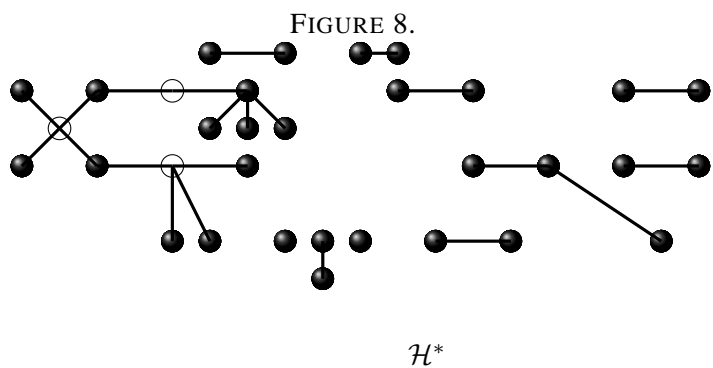
\begin{figure}[h] 
\caption{}\label{EstarReduce}
\begin{center}
\begin{tikzpicture}
\shade [shading=ball, ball color=black] (-2,1) circle (.15);
\shade [shading=ball, ball color=black] (-2,0) circle (.15);
\shade [shading=ball, ball color=black] (0.5,-1) circle (.15);
\shade [shading=ball, ball color=black] (0,-1) circle (.15);
\shade [shading=ball, ball color=black] (-1,1) circle (.15);
\shade [shading=ball, ball color=black] (-1,0) circle (.15);
\draw  [shape=circle] (-1.5,0.5) circle (.15);
\draw  [shape=circle] (0,0) circle (.15);
\draw  [shape=circle] (0,1) circle (.15);
\shade [shading=ball, ball color=black] (1,0) circle (.15);

\shade [shading=ball, ball color=black] (4,0) circle (.15);
\shade [shading=ball, ball color=black]  (5,0) circle (.15);
\shade [shading=ball, ball color=black] (6,0) circle (.15);
\shade [shading=ball, ball color=black]  (7,0) circle (.15);
\shade [shading=ball, ball color=black] (0.5,0.5) circle (.15);
\shade [shading=ball, ball color=black] (1.5,0.5) circle (.15);
\shade [shading=ball, ball color=black] (1,0.5) circle (.15);
\shade [shading=ball, ball color=black] (1.5,-1) circle (.15);
\shade [shading=ball, ball color=black] (2,-1) circle (.15);
\shade [shading=ball, ball color=black] (2,-1.5) circle (.15);

\shade [shading=ball, ball color=black]  (2.5,-1) circle (.15);
\shade [shading=ball, ball color=black] (3.5,-1) circle (.15);
\shade [shading=ball, ball color=black]  (4.5,-1) circle (.15);
\shade [shading=ball, ball color=black] (1,1) circle (.15);

\shade [shading=ball, ball color=black](3,1) circle (.15);
\shade [shading=ball, ball color=black] (4,1) circle (.15);

\shade [shading=ball, ball color=black] (6,1) circle (.15);
\shade [shading=ball, ball color=black]  (7,1) circle (.15);

\shade [shading=ball, ball color=black] (1.5,1.5) circle (.15);
\shade [shading=ball, ball color=black] (0.5,1.5) circle (.15);

\shade [shading=ball, ball color=black] (2.5,1.5) circle (.15);
\shade [shading=ball, ball color=black]  (3,1.5) circle (.15);
\shade [shading=ball, ball color=black]  (6.5,-1) circle (.15);

\draw [line width=1.2pt ] (1,1)--(.5,0.5)  ;
\draw [line width=1.2pt ] (1,1)--(1.5,0.5)  ;
\draw [line width=1.2pt ] (1,1)--(1,0.5)  ;
\draw [line width=1.2pt ] (-1.5,0.5)--(-2,0)  ;
\draw [line width=1.2pt ] (-1.5,0.5)--(-2,1)  ;
\draw [line width=1.2pt ] (0,0)--(0.5,-1)  ;
\draw [line width=1.2pt ] (0,0)--(0,-1)  ;
\draw [line width=1.2pt ] (-1.5,0.5)--(-1,0)  ;
\draw [line width=1.2pt ] (-1.5,0.5)--(-1,1)  ;
\draw [line width=1.2pt ] (0,0)--(-1,0)  ;
\draw [line width=1.2pt ] (0,1)--(-1,1)  ;
\draw [line width=1.2pt ] (0,1)--(1,1)  ;
\draw [line width=1.2pt ] (0,0)--(1,0)  ;

\draw [line width=1.2pt ] (4,0)--(5,0)  ;
\draw [line width=1.2pt ] (6,0)--(7,0)  ;

\draw [line width=1.2pt ] (7,1)--(6,1)  ;

\draw [line width=1.2pt ] (3,1)--(4,1)  ;

\draw [line width=1.2pt ] (0.5,1.5)--(1.5,1.5)  ;

\draw [line width=1.2pt ] (3,1.5)--(2.5,1.5)  ;

\draw [line width=1.2pt ] (2,-1)--(2,-1.5)  ;

\draw [line width=1.2pt ] (3.5,-1)--(4.5,-1)  ;

\draw [line width=1.2pt] (6.5,-1)--(5,0)  ;
\path (3,-2.5)--(3,-2.5) node [pos=.5, right ] {$\H^*$};

\end{tikzpicture}
\end{center}
\end{figure}

We now prove an analogous formula for the projective dimension of $\H$ when the 1-dimensional part, $\H^1$, of $\H$ is a string or cycle of $1$-stars. Notice that the hypergraph needs not be 1-dimensional itself, it suffices that its 1-skeleton satisfies certain properties and the higher dimensional faces are ``well behaved''. 

\begin{Proposition}\label{StarString} Let $\H$ be a hypergraph with $V(\H)=V(\H^1)$; assume its 1-skeleton $\H^1$ is a separated string or cycle of $1$-stars, and distinct stars do not share any higher dimensional face.
\begin{enumerate}
\item If the 1-skeleton of $\H$ is an open cycle of $n_1$ 1-stars, then $\pd(\H)=|V(\H)|-1-\left\lfloor \frac{n_1-2}{3}\right\rfloor$;
\item in all other cases $\pd(\H)=|V(\H)|-s^*-\sum_{i=1}^{s^*}\left\lfloor \frac{n_{i}-1}{3}\right\rfloor +M^*(\H)$,
where $s^*$ is the number of strings of open stars in $\H^1$, $n_{i}$ is the number of open stars in each string of open stars.
\end{enumerate} 
\end{Proposition}

\demo
By Proposition \ref{eq}, we may assume that there are no adjacent closed vertices. Moreover, by definition, all 1-stars are open stars. We induct on the number of 1-stars and number of vertices. If there is no 1-stars, then $\H^1$ is either a string or a cycle hypergraph; the assumption on the higher dimensional faces implies that $\H=\H^1$ (the distinct vertices of $\H$ are $0$-stars); also, since $\H^1$ is separated, the formula follows by \cite[Theorem~3.4, 4.3]{LMa}.

 We may then assume there is at least one proper $1$-star; let $w$ be its open center and $v$ an endpoint of one of its branches. Let $F_1,\ldots,F_t$ be all higher dimensional faces containing $v$, we apply Lemma \ref{ref} to $v$ and observe that $\H_{v}=\Q_{v}\cup\{w\}\cup(\bigcup_i\{F_i|_{\H_v}\})$, because $w$ is a joint and by the assumption on the higher dimensional faces. Then all the other vertices of each $F_i|_{\H_v}$ are closed in $\H_v$, hence by Proposition \ref{eq}, we have $\pd(\H_{v})=\pd(\Q_{v}\cup\{w\})$.  Since $w$ is open in $\Q_{v}$, then by Lemma \ref{H^0}, we have $\pd(\H_{v})\leq\pd(\Q_{v})+1$. This inequality combined with Lemma \ref{ref} yields $\pd(\H)=\pd(\Q_v)+1$. 
Then by induction (and since $Q_{v}$ has $|V(\H)|-1$ vertices)
$$\pd(\H)=\pd(\Q_{v})+1=(|V(\H)|-1)-s^*-\sum_{i=1}^{s^*}\left\lfloor \frac{n_{i}-1}{3}\right\rfloor +M^*(\H)+1,$$
yielding the desired formula.
\QED
\bigskip

\section{Appendix: algorithmic procedures and more examples}

In order to apply Propositions \ref{AttachString} and \ref{cutTree} to compute the projective dimension of a hypergraph, we first need to algorithmically recognize if a vertex $i$ in $\H$ is a joint or an endpoint, thus we need to determine its degree $d(i)$. Actually, for Algorithm \ref{AlgPd} it suffices to know if $d(i)=0$, $1$, $2$ or if it is greater than $2$; so, for reasons of efficiency (e.g. if $i$ has a large number of neighbors), in Algorithm \ref{deg} we only consider these possible outputs; of course, it can be easily modified to actually compute $d(i)$. The auxiliary variable $j$ runs through the elements of the vertex set to identify neighbors of $i$.

\begin{Algorithm}\label{deg}
Let $\H$ be a hypergraph, $V(\mathcal{H})=\{1,2,\cdots,\mu\}$.
The input is: $i\in V(\mathcal{H})$, i.e. a vertex in a hypergraph $\mathcal{H}$. The output is: $n=d(i)$, if this number is $0,1,2$, or ``$n>2$'' otherwise.

\begin{itemize}
\item[Step 0:] Set $n=0$, $V=V(\mathcal{H})$ and $j=1$.

\item[Step 1:] If $n=3$, then stop and give ``$n>2$'' as output.\\
If $|V|=1$ then stop and give $n$ as output;\\
If $j=\mu$, then stop and give $n$ as output.\\
Otherwise, go to Step 2.

\item[Step 2:] If $j=i$, set $j=j+1$ and go to Step 1.\\
 If $j\neq i$ then set $V=V\backslash\{j\}$ and do the following: if $\{i,j\}\in\mathcal{H}$, then set $n=n+1$ and go to Step 1.  If $\{i,j\}\notin\mathcal{H}$ go to Step 1.
\end{itemize}
\end{Algorithm}

The following result provides an effective algorithmic way to compute the projective dimension of each connected component in Theorem \ref{uniqueCycle}. In the following algorithm we use the variable $i$ to detect the vertices with degree one (if any); the variable $j$ runs through the other vertices looking for neighbors of $i$, and $k$ looks for the other neighbor of $j$ (if any). The variable $v$ is used to count $|V(\H)|$ (as the algorithm runs $\H$ changes and so does $|V(\H)|$), and $c$ is used to isolate the scenario where $\H$ is a $v$-cycle.

\begin{Algorithm}\label{AlgPd}
Input: A connected 1-dimensional hypergraph $\mathcal{H}$ with
at most one cycle. Let the vertex set be $V(\mathcal{H})=\{1,2,\cdots,\mu\}$. The output is: $P=\pd(\mathcal{H})$.

\begin{itemize}

\item [Step 0:] Set $P=0$, $v=\mu$ and $i=1$.

\item [Step 1:] If $\mathcal{H}=\emptyset$, stop the process and give $P$ as output.\\
If $\H\neq \emptyset$ set $j=k=1$, $c=0$ and do the following: if $i\leq \mu$,  and go to Step 2, if $i=\mu+1$, then set $i=1$ and go to Step 2.
\item [Step 2:] If $i\notin V(\H)$, then set $i=i+1$ and start Step 2 again.\\
If $i\in V(\H)$, compute $d(i)$ using Algorithm \ref{deg}.\\
\indent if $d(i)=0$, set $\mathcal{H}=\mathcal{H}_{i}$,  $P=P+1$, $v=v-1$ and $i=i+1$, then go to Step 1;\\
\indent if $d(i)=1$ then go to Step 3;\\
\indent if $d(i)>1$ set $c=c+1$ and do the following:\\
\indent\indent\indent if $c=v$, then $\H$ is a $v$-cycle and we go to Step 7;\\
\indent\indent\indent if $c<v$ then do the following:\\
\indent\indent\indent\indent\indent if $i\neq \mu$, set $i=i+1$ and start Step 2 again;\\
\indent\indent\indent\indent\indent if $d(i)>1$ and $i=\mu$, set $i=1$ and start Step 2 again.

\item [Step 3:] If $j=i$ or if $\{i,j\}\notin\mathcal{H}$ then set $j=j+1$ and start again Step 3. If $\{i,j\}\in\mathcal{H}$ then go to Step 4. 

\item [Step 4:] Check if $\{j\}\in\mathcal{H}$. If so, set $\mathcal{H}=\mathcal{H}_{i}$ and $P=P+1$, $v=v-1$, $i=i+1$ then go to Step 1. if If $\{j\}\notin\mathcal{H}$ go to Step 5. (notice that since $\{j\}\notin \mathcal{H}$, then $d(j)\geq 2$.)

\item [Step 5:] Use Algorithm \ref{deg} to compute $d(j)$. 
 If $d(j)=2$, then go to Step 6; otherwise set $\mathcal{H}=\mathcal{H}\backslash\{\{i,j\},\{i\}\}$, $P=P+1$, $v=v-1$, $i=i+1$ and go to Step 1.

\item [Step 6:] If $k=i$, or if $k=j$, or if $\{j,k\}\notin\mathcal{H}$, then set $k=k+1$ and start again Step 6. If $\{j,k\}\in\mathcal{H}$ then set $\mathcal{H}=\mathcal{H}_{i,j,k}$, $P=P+2$, $v=v-3$, $i=i+1$ and go to Step 1. (this procedure stops because $d(j)=2$)

\item [Step 7:] Use Algorithm 5.6 in \cite{LMa} to compute $\pd(\mathcal{H})=P_{c}$.
Output $P=P+P_{c}$. 
\end{itemize}
\end{Algorithm}

\begin{Remark}
$($1$)$ The variable $c$ counts the number of times that Step 2 runs consecutively without finding a vertex with degree $\leq 1$. If $c=|V(\H)|$, then every vertex of $\H$ has degree $\geq 2$, so $\H$ is a cycle. \\
$($2$)$ Step 3 always starts with $j=1$, and since $d(i)=1$ then there is precisely one $j$ with $\{i,j\}\in \H$; therefore Step 3 does not need a line for the case where $j$ becomes larger than $\mu$, because it stops before then. A similar comment holds for the variable $k$ in Step 6. \\
$($3$)$ In Step 4 of Algorithm \ref{AlgPd}, we set $\H=\H_i$ because $i$ has only one neighbor (so $i$ is closed in $\H$), which is also closed. Thus, by Lemma \ref{red}, we can remove the vertex $i$ and add one to $P$.
\end{Remark}

The following example illustrates the use of Algorithm \ref{AlgPd}.

\begin{Example}: In Figure \ref{algEx1} we provide a hypergraph $\mathcal{H}$ and all the steps of Algorithm \ref{AlgPd} to compute $\pd(\H)$.

\begin{figure}[h] 
\caption{}\label{algEx1}
\begin{center}
\begin{tikzpicture}[thick, scale=0.65]
\tikzstyle{ghostfill} = [fill=white]
         \tikzstyle{ghostdraw} = [draw=black!50]
\usetikzlibrary{arrows,shapes,positioning}
\usetikzlibrary{decorations.markings}
\tikzstyle arrowstyle=[scale=1]
\tikzstyle directed=[postaction={decorate,decoration={markings,
    mark=at position .65 with {\arrow[arrowstyle]{stealth}}}}]
\tikzstyle reverse directed=[postaction={decorate,decoration={markings,
    mark=at position .65 with {\arrowreversed[arrowstyle]{stealth};}}}]

\shade [shading=ball, ball color=black] (0,0) circle (.15);
\draw  [shape=circle] (1,0) circle (.15);
\draw  [shape=circle] (2,0) circle (.15);
\draw  [shape=circle] (3,0) circle (.15);
\draw  [shape=circle] (4,0) circle (.15);
\draw  [shape=circle] (5,0) circle (.15);
\draw  [shape=circle] (6,0) circle (.15);
\draw  [shape=circle] (7,0) circle (.15);
\shade [shading=ball, ball color=black] (8,0) circle (.15);
\shade [shading=ball, ball color=black] (3,1) circle (.15);
\draw  [shape=circle] (4,1) circle (.15);
\draw  [shape=circle] (5,1) circle (.15);
\draw  [shape=circle] (6,1) circle (.15);
\draw  [shape=circle] (7,1) circle (.15);
\shade [shading=ball, ball color=black] (8,1) circle (.15) node [ above] {$j$};
\shade [shading=ball, ball color=black] (9,1) circle (.15)node [ above] {$i$};
\shade [shading=ball, ball color=black] (2,0.5) circle (.15);

\shade [shading=ball, ball color=black] (6,0.5) circle (.15);
\shade [shading=ball, ball color=black]  (6,-0.5) circle (.15);

\draw [line width=1.2pt ] (0,0)--(1,0)  ;
\draw [line width=1.2pt ] (1,0)--(2,0)  ;
\draw [line width=1.2pt ] (2,0)--(3,0)  ;
\draw [line width=1.2pt ] (3,0)--(4,0)  ;
\draw [line width=1.2pt ] (4,0)--(5,0)  ;
\draw [line width=1.2pt ] (5,0)--(6,0)  ;
\draw [line width=1.2pt ] (6,0)--(7,0)  ;
\draw [line width=1.2pt ] (7,0)--(8,0)  ;
\draw [line width=1.2pt ] (3,1)--(4,1)  ;
\draw [line width=1.2pt ] (4,1)--(5,1)  ;
\draw [line width=1.2pt ] (5,1)--(6,1)  ;
\draw [line width=1.2pt ] (6,1)--(7,1)  ;
\draw [line width=1.2pt ] (7,1)--(8,1)  ;
\draw [line width=1.2pt ] (8,1)--(9,1)  ;

\draw [line width=1.2pt ] (3,0)--(3,1)  ;
\draw [line width=1.2pt ] (4,0)--(4,1)  ;
\draw [line width=1.2pt ] (2,0)--(2,0.5)  ;
\draw [line width=1.2pt ] (6,0)--(6,0.5)  ;
\draw [line width=1.2pt ] (6,0)--(6,-0.5)  ;

\path (4,-1)--(4,-1) node [pos=.5, right ] {$\H$, $P=0$};
 \draw[thick, ->, color=red] (10,0.5)--(11,0.5);

\shade [shading=ball, ball color=black] (12,0) circle (.15);
\draw  [shape=circle] (13,0) circle (.15);
\draw  [shape=circle] (14,0) circle (.15);
\draw  [shape=circle] (15,0) circle (.15);
\draw  [shape=circle] (16,0) circle (.15);
\draw  [shape=circle] (17,0) circle (.15);
\draw  [shape=circle] (18,0) circle (.15);
\draw  [shape=circle] (19,0) circle (.15);
\shade [shading=ball, ball color=black] (20,0) circle (.15);
\shade [shading=ball, ball color=black] (15,1) circle (.15);
\draw  [shape=circle] (16,1) circle (.15);
\draw  [shape=circle] (17,1) circle (.15);
\draw  [shape=circle] (18,1) circle (.15) node [ above] {$k$};
\draw  [shape=circle] (19,1) circle (.15) node [ above] {$j$};
\shade [shading=ball, ball color=black] (20,1) circle (.15) node [ above] {$i$};

\shade [shading=ball, ball color=black] (14,.5) circle (.15);
\shade [shading=ball, ball color=black] (18,.5) circle (.15);
\shade [shading=ball, ball color=black]  (18,-.5) circle (.15);

\draw [line width=1.2pt ] (12,0)--(13,0)  ;
\draw [line width=1.2pt ] (13,0)--(14,0)  ;
\draw [line width=1.2pt ] (14,0)--(15,0)  ;
\draw [line width=1.2pt ] (15,0)--(16,0)  ;
\draw [line width=1.2pt ] (16,0)--(17,0)  ;
\draw [line width=1.2pt ] (17,0)--(18,0)  ;
\draw [line width=1.2pt ] (18,0)--(19,0)  ;
\draw [line width=1.2pt ] (19,0)--(20,0)  ;
\draw [line width=1.2pt ] (15,1)--(16,1)  ;
\draw [line width=1.2pt ] (16,1)--(17,1)  ;
\draw [line width=1.2pt ] (17,1)--(18,1)  ;
\draw [line width=1.2pt ] (18,1)--(19,1)  ;
\draw [line width=1.2pt ] (19,1)--(20,1)  ;

\draw [line width=1.2pt ] (15,0)--(15,1)  ;
\draw [line width=1.2pt ] (16,0)--(16,1)  ;
\draw [line width=1.2pt ] (14,0)--(14,.5)  ;

\draw [line width=1.2pt ] (18,0)--(18,0.5)  ;
\draw [line width=1.2pt ] (18,0)--(18,-.5)  ;

\path (15,-1)--(15,-1) node [pos=.5, right ] {$\H=\H_i$, $P=1$};
\draw[thick, ->, color=red] (15,-1)--(15,-2);

\shade [shading=ball, ball color=black] (12,-3.5) circle (.15);
\draw  [shape=circle] (13,-3.5) circle (.15);
\draw  [shape=circle] (14,-3.5) circle (.15);
\draw  [shape=circle] (15,-3.5) circle (.15);
\draw  [shape=circle] (16,-3.5) circle (.15);
\draw  [shape=circle] (17,-3.5) circle (.15);
\draw  [shape=circle] (18,-3.5) circle (.15);
\draw  [shape=circle] (19,-3.5) circle (.15);
\shade [shading=ball, ball color=black] (20,-3.5) circle (.15);
\shade [shading=ball, ball color=black] (15,-2.5) circle (.15);
\draw  [shape=circle] (16,-2.5) circle (.15);
\shade [shading=ball, ball color=black] (17,-2.5) circle (.15) node [above] {$i$};

\shade [shading=ball, ball color=black] (14,-3) circle (.15);
\shade [shading=ball, ball color=black] (18,-3) circle (.15);
\shade [shading=ball, ball color=black]  (18,-4) circle (.15);

\draw [line width=1.2pt ] (12,-3.5)--(13,-3.5)  ;
\draw [line width=1.2pt ] (13,-3.5)--(14,-3.5)  ;
\draw [line width=1.2pt ] (14,-3.5)--(15,-3.5)  ;
\draw [line width=1.2pt ] (15,-3.5)--(16,-3.5)  ;
\draw [line width=1.2pt ] (16,-3.5)--(17,-3.5)  ;
\draw [line width=1.2pt ] (17,-3.5)--(18,-3.5)  ;
\draw [line width=1.2pt ] (18,-3.5)--(19,-3.5)  ;
\draw [line width=1.2pt ] (19,-3.5)--(20,-3.5)  ;
\draw [line width=1.2pt ] (15,-2.5)--(16,-2.5)  ;
\draw [line width=1.2pt ] (16,-2.5)--(17,-2.5)  ;
\draw [line width=1.2pt ] (15,-3.5)--(15,-2.5)  ;
\draw [line width=1.2pt ] (16,-3.5)--(16,-2.5)  ;
\draw [line width=1.2pt ] (14,-3.5)--(14,-3)  ;
\draw [line width=1.2pt ] (18,-3.5)--(18,-3)  ;
\draw [line width=1.2pt ] (18,-3.5)--(18,-4)  ;

\path (15,-4.5)--(15,-4.5) node [pos=.5, right ] {$\H=\H_{i,j,k}$, $P=1+2=3$};

\draw[thick, <-, color=red] (10,-3)--(11,-3);

\shade [shading=ball, ball color=black] (0,-3.5) circle (.15) node[below]{$i$};
\draw  [shape=circle] (1,-3.5) circle (.15)node[below]{$j$};
\draw  [shape=circle] (2,-3.5) circle (.15)node[below]{$k$};
\draw  [shape=circle] (3,-3.5) circle (.15);
\draw  [shape=circle] (4,-3.5) circle (.15);
\draw  [shape=circle] (5,-3.5) circle (.15);
\draw  [shape=circle] (6,-3.5) circle (.15);
\draw  [shape=circle] (7,-3.5) circle (.15);
\shade [shading=ball, ball color=black] (8,-3.5) circle (.15);
\shade [shading=ball, ball color=black] (3,-2.5) circle (.15);
\draw  [shape=circle] (4,-2.5) circle (.15);

\shade [shading=ball, ball color=black] (2,-3) circle (.15);
\shade [shading=ball, ball color=black] (6,-3) circle (.15);
\shade [shading=ball, ball color=black]  (6,-4) circle (.15);

\draw [line width=1.2pt ] (0,-3.5)--(1,-3.5)  ;
\draw [line width=1.2pt ] (1,-3.5)--(2,-3.5)  ;
\draw [line width=1.2pt ] (2,-3.5)--(3,-3.5)  ;
\draw [line width=1.2pt ] (3,-3.5)--(4,-3.5)  ;
\draw [line width=1.2pt ] (4,-3.5)--(5,-3.5)  ;
\draw [line width=1.2pt ] (5,-3.5)--(6,-3.5)  ;
\draw [line width=1.2pt ] (6,-3.5)--(7,-3.5)  ;
\draw [line width=1.2pt ] (7,-3.5)--(8,-3.5)  ;
\draw [line width=1.2pt ] (3,-2.5)--(4,-2.5)  ;
\draw [line width=1.2pt ] (3,-3.5)--(3,-2.5)  ;
\draw [line width=1.2pt ] (4,-3.5)--(4,-2.5)  ;
\draw [line width=1.2pt ] (2,-3.5)--(2,-3)  ;
\draw [line width=1.2pt ] (6,-3.5)--(6,-3)  ;
\draw [line width=1.2pt ] (6,-3.5)--(6,-4)  ;

\path (2.5,-4.5)--(2.5,-4.5) node [pos=.5, right ] {$\H=\mathcal{H}\backslash\{\{i,j\},\{i\}\}$};
\path (2.5,-5.2)--(2.5,-5.2) node [pos=.5, right ] {$P=3+1=4$};
\draw[thick, <-, color=red] (2.5,-5.5)--(2.5,-4.5);

\shade [shading=ball, ball color=black] (3,-7) circle (.15);
\draw  [shape=circle] (4,-7) circle (.15);
\draw  [shape=circle] (5,-7) circle (.15);
\draw  [shape=circle] (6,-7) circle (.15);
\draw  [shape=circle] (7,-7) circle (.15);
\shade [shading=ball, ball color=black] (8,-7) circle (.15);
\shade [shading=ball, ball color=black] (3,-6) circle (.15);
\draw  [shape=circle] (4,-6) circle (.15);
\shade [shading=ball, ball color=black] (2,-6.5) circle (.15) node [above]{$i$};
\shade [shading=ball, ball color=black] (6,-6.5) circle (.15);
\shade [shading=ball, ball color=black]  (6,-7.5) circle (.15);

\draw [line width=1.2pt ] (3,-7)--(4,-7)  ;
\draw [line width=1.2pt ] (4,-7)--(5,-7)  ;
\draw [line width=1.2pt ] (5,-7)--(6,-7)  ;
\draw [line width=1.2pt ] (6,-7)--(7,-7)  ;
\draw [line width=1.2pt ] (7,-7)--(8,-7)  ;
\draw [line width=1.2pt ] (3,-6)--(4,-6)  ;
\draw [line width=1.2pt ] (3,-7)--(3,-6)  ;
\draw [line width=1.2pt ] (4,-7)--(4,-6)  ;

\draw [line width=1.2pt ] (6,-7)--(6,-6.5)  ;
\draw [line width=1.2pt ] (6,-7)--(6,-7.5)  ;

\path (2.5,-8)--(2.5,-8) node [pos=.5, right ] {$\H=\H_{i,j,k}$, $P=4+2=5$};
\draw[thick, ->, color=red] (10,-6.5)--(11,-6.5);

\shade [shading=ball, ball color=black] (15,-7) circle (.15);
\draw  [shape=circle] (16,-7) circle (.15);
\draw  [shape=circle] (17,-7) circle (.15);
\draw  [shape=circle] (18,-7) circle (.15) node[above right]{$k$};
\draw  [shape=circle] (19,-7) circle (.15)node[above right]{$j$};
\shade [shading=ball, ball color=black] (20,-7) circle (.15)node[above right]{$i$};
\shade [shading=ball, ball color=black] (15,-6) circle (.15);
\draw  [shape=circle] (16,-6) circle (.15);
\shade [shading=ball, ball color=black] (18,-6.5) circle (.15);
\shade [shading=ball, ball color=black]  (18,-7.5) circle (.15);

\draw [line width=1.2pt ] (15,-7)--(16,-7)  ;
\draw [line width=1.2pt ] (16,-7)--(17,-7)  ;
\draw [line width=1.2pt ] (17,-7)--(18,-7)  ;
\draw [line width=1.2pt ] (18,-7)--(19,-7)  ;
\draw [line width=1.2pt ] (19,-7)--(20,-7)  ;
\draw [line width=1.2pt ] (15,-6)--(16,-6)  ;
\draw [line width=1.2pt ] (15,-7)--(15,-6)  ;
\draw [line width=1.2pt ] (16,-7)--(16,-6)  ;
\draw [line width=1.2pt ] (18,-7)--(18,-6.5)  ;
\draw [line width=1.2pt ] (18,-7)--(18,-7.5)  ;

\path (15,-8)--(15,-8) node [pos=.5, right ] {$\H=\H_{i}$, $P=5+1=6$};
\draw[thick, ->, color=red] (15,-8)--(15,-9);

\shade [shading=ball, ball color=black] (15,-10.5) circle (.15);
\draw  [shape=circle] (16,-10.5) circle (.15);
\shade [shading=ball, ball color=black] (17,-10.5) circle (.15);
\shade [shading=ball, ball color=black] (15,-9.5) circle (.15);
\draw  [shape=circle] (16,-9.5) circle (.15);
\shade [shading=ball, ball color=black] (18,-10) circle (.15);
\shade [shading=ball, ball color=black]  (18,-11) circle (.15);
\draw [line width=1.2pt ] (15,-10.5)--(16,-10.5)  ;
\draw [line width=1.2pt ] (16,-10.5)--(17,-10.5)  ;
\draw [line width=1.2pt ] (15,-9.5)--(16,-9.5)  ;
\draw [line width=1.2pt ] (15,-10.5)--(15,-9.5)  ;
\draw [line width=1.2pt ] (16,-10.5)--(16,-9.5)  ;

\path (15,-11.5)--(15,-11.5) node [pos=.5, right ] {$\H=\H_{i,j,k}$, $P=6+2=8$};
\draw[thick, <-, color=red] (10,-10)--(11,-10) node [right]{...};

\shade [shading=ball, ball color=black] (3,-10.5) circle (.15);
\draw  [shape=circle] (4,-10.5) circle (.15);
\shade [shading=ball, ball color=black] (3,-9.5) circle (.15);
\draw  [shape=circle] (4,-9.5) circle (.15);
\draw [line width=1.2pt ] (3,-10.5)--(4,-10.5)  ;
\draw [line width=1.2pt ] (3,-9.5)--(4,-9.5)  ;
\draw [line width=1.2pt ] (3,-10.5)--(3,-9.5)  ;
\draw [line width=1.2pt ] (4,-10.5)--(4,-9.5)  ;

\path (2,-11.5)--(2,-11.5) node [pos=.5, right ] {$P=P_c+11=3+11=14$};

\end{tikzpicture}
\end{center}
\end{figure}

\end{Example}


\begin{thebibliography}{99}

\bibitem{Be}{ C. Berge, Hypergraphs: Combinatorics of Finite Sets. North–Holland, Amsterdam, 1989.}

\bibitem{C}{ D. Cook II, The uniform face ideals of a simplicial complex, J. Commut. Algebra (to appear), 52
pages; arXiv:1308.1299.}

\bibitem{DHS}{ H. Dao, C. Huneke, J. Schweig, Bounds on the regularity and projective dimension of ideals associated to graphs, J. Algebraic Combin. {\bf 38} (2013), no. 1, 37-55.}


\bibitem{DS}{ H. Dao, J. Schweig, Projective dimension, graph domination parameters, and independence complex
homology, J. Combin. Theory Ser. A, {\bf 120} (2013), no. 2,  453-469}

\bibitem{DS1}{ H. Dao, J. Schweig, Further applications of clutter domination parameters to projective dimension (with H. Dao), J. Algebra {\bf 432} (2015), 1-11.}


\bibitem{F}{ S. Faridi, The projective dimension of sequentially Cohen-Macaulay monomial ideals, arXiv:1310.5598 (2013)}


\bibitem{H}{ H. T. H\`a, Regularity of squarefree monomial ideals. In: Cooper, S.M., Sather-Wagstaff, S. (eds.)
Connections Between Algebra, Combinatorics, and Geometry. Springer Proceedings in Mathematics \& Statistics, vol. {\bf 76} (2014), 251-276}


\bibitem{HL}{ H. T. H\`a and K.-N. Lin, Normal 0-1 polytopes, SIAM J. Discrete Math., {\bf 29} (2015), no. 1, 210-223.}

\bibitem{HW}{ H. T. H\`a and R. Woodroofe, Results on the regularity of square-free monomial ideals. Adv. in Appl. Math. {\bf 58} (2014), 21-36.}

\bibitem{KRT}{ K. Kimura, G. Rinaldo, N. Terai, Arithmetical rank of squarefree monomial ideals generated by five elements or with arithmetic degree four, Comm. Algebra. {\bf 40} (2012), 4147-4170.}


\bibitem{KTY}{ K. Kimura, N. Terai and K. Yoshida, Arithmetical rank of square-free monomial ideals
of small arithmetic degree, J. Algebra Comb. {\bf 29} (2009), 389-404.}

\bibitem{KTY1}{ K. Kimura, N. Terai and K. Yoshida, Arithmetical rank of a squarefree monomial ideal whose Alexander dual is of deviation two. (English summary) 
Acta Math. Vietnam. {\bf 40} (2015), no. 3, 375-391.} 

\bibitem{KM}{ K. Kimura and P. Mantero, Arithmetical Rank of strings and cycles, to appear in J. Commut. Algebra, arXiv:1407.5571}

\bibitem{Ku}{ M. Kummini, Regularity, depth and arithmetic rank of bipartite edge ideals, J. Algebraic Combin. 30, no. {\bf 4} (2009), 429-445.}

\bibitem{LMa}{ K.-N. Lin and P. Mantero, Projective Dimension of Strings and Cycles, to appear in Comm. in Algebra.}

\bibitem{LMc}{ K.-N. Lin and J. McCullough, Hypergraphs and regularity of square-free monomial ideals, Internat. J. Algebra Comput. {\bf 23} (2013), 1573-1590.}

\bibitem{M2} D.R. Grayson, M.E. Stillman, Macaulay 2, a software system for research in algebraic geometry.
\newblock \verb|http://www.math.uiuc.edu/Macaulay2/|.

\bibitem{MV}{ S. Morey and R.H. Villarreal, Edge ideals: algebraic and combinatorial properties. Progress in commutative algebra 1, 85-126, de Gruyter, Berlin, 2012.}


\bibitem{T}{D. Taylor,  Ideals generated by monomials in an $R$--sequence, Ph.\ D.\ thesis, University of Chicago, 1966.}
\end{thebibliography}
\end{document}